\documentclass[11pt]{amsart}
\usepackage{amscd,amssymb}
\usepackage[all]{xy}
\usepackage[colorlinks,plainpages,backref,urlcolor=blue]{hyperref}

\newtheorem{theorem}{Theorem}[section]
\newtheorem{corollary}[theorem]{Corollary}
\newtheorem{lemma}[theorem]{Lemma}
\newtheorem{prop}[theorem]{Proposition}

\newtheorem{claim}{Claim}

\theoremstyle{definition}

\newtheorem{example}[theorem]{Example}
\newtheorem{remark}[theorem]{Remark}
\newtheorem*{ack}{Acknowledgments}

\newenvironment{romenum}
{

\begin{enumerate}}{\end{enumerate}}

\newenvironment{alphenum}
{

\begin{enumerate}}{\end{enumerate}}

\newcommand{\Z}{\mathbb{Z}}
\newcommand{\Q}{\mathbb{Q}}
\newcommand{\R}{\mathbb{R}}

\newcommand{\KK}{\mathbb{K}}
\newcommand{\RP}{\mathbb{RP}}
\newcommand{\CP}{\mathbb{CP}}
\newcommand{\LL}{\mathbb{L}}
\renewcommand{\k}{\Bbbk}

\newcommand{\h}{\mathfrak{h}}
\newcommand{\g}{\mathfrak{g}}
\newcommand{\M}{\mathfrak{m}}

\newcommand{\K}{\mathfrak{K}}

\newcommand{\RR}{\mathcal{R}}
\newcommand{\VV}{\mathcal{V}}

\newcommand{\PP}{\mathcal{P}}

\newcommand{\E}{\mathsf{E}}
\newcommand{\V}{\mathsf{V}}
\newcommand{\sE}{\mathsf{E}}
\newcommand{\sV}{\mathsf{V}}
\newcommand{\sW}{\mathsf{W}}

\DeclareMathOperator{\rank}{rank}
\DeclareMathOperator{\gr}{gr}
\DeclareMathOperator{\im}{im}
\DeclareMathOperator{\coker}{coker}

\DeclareMathOperator{\id}{id}
\DeclareMathOperator{\ab}{{ab}}

\DeclareMathOperator{\ch}{char}

\DeclareMathOperator{\Hom}{{Hom}}
\DeclareMathOperator{\Tor}{{Tor}}

\DeclareMathOperator{\spn}{span}

\DeclareMathOperator{\ev}{ev}

\DeclareMathOperator{\ord}{{ord}}
\DeclareMathOperator{\supp}{{supp}}
\DeclareMathOperator{\lk}{lk}
\DeclareMathOperator{\Tors}{{Tors}}
\DeclareMathOperator{\Lie}{{Lie}}
\DeclareMathOperator{\diag}{{diag}}
\DeclareMathOperator{\FP}{{FP}}

\newcommand{\wL}{\widehat{\Lambda}}
\newcommand{\wh}{\widehat{\h}}
\renewcommand{\L}{\Lambda}
\newcommand{\G}{\Gamma}

\newcommand{\same}{\Leftrightarrow}
\newcommand{\surj}{\twoheadrightarrow}
\newcommand{\inj}{\hookrightarrow}
\newcommand{\isom}{\xrightarrow{\,\simeq\,}}
\def\set#1{{\left\{#1\right\}}}
\newcommand{\abs}[1]{\left| #1 \right|}

\newcommand{\SR}{\k\langle L\rangle}
\newcommand{\chik}{\chi^{\:}_{\k}}

\begin{document}

\title[Toric complexes and Artin kernels]{%
Toric complexes and Artin kernels}

\author[Stefan Papadima]{Stefan Papadima$^1$}
\address{Institute of Mathematics Simion Stoilow, 
P.O. Box 1-764,
RO-014700 Bucharest, Romania}
\email{Stefan.Papadima@imar.ro}
\thanks{$^1$Partially supported by the CEEX Programme of the 
Romanian Ministry of Education and Research, contract
2-CEx 06-11-20/2006}

\author[Alexander~I.~Suciu]{Alexander~I.~Suciu}
\address{Department of Mathematics,
Northeastern University,
Boston, MA 02115, USA}
\email{a.suciu@neu.edu}



\begin{abstract}
A simplicial complex $L$ on $n$ vertices determines 
a subcomplex $T_L$ of the $n$-torus, with fundamental 
group the right-angled Artin group $G_{L}$.  Given an 
epimorphism $\chi\colon G_{L}\to \Z$, let $T_L^{\chi}$ 
be the corresponding cover, with fundamental group 
the Artin kernel $N_{\chi}$. We compute the cohomology 
jumping loci of the toric complex $T_L$, as well as the 
homology groups of $T_L^{\chi}$ with coefficients in a 
field $\k$, viewed as modules over the group algebra $\k\Z$.  
We give combinatorial conditions for $H_{\le r}(T_L^{\chi};\k)$ 
to have trivial $\Z$-action, allowing us to compute  
the truncated cohomology ring, $H^{\le r}(T_L^{\chi};\k)$.  
We also determine several Lie algebras associated to 
Artin kernels, under certain triviality assumptions on the 
monodromy $\Z$-action, and establish the $1$-formality 
of these (not necessarily finitely presentable) groups. 
\end{abstract}
\maketitle

\tableofcontents

\newpage
\section{Introduction}
\label{sect:intro}

The underlying theme of this paper is the interplay between 
topology, discrete geometry, group theory, and commutative 
algebra, as revealed by the intricate connections tying up the 
combinatorics of a simplicial complex with the algebraic and 
geometric topology of certain spaces modeled on it. 

Our first goal is to better understand the topology of a toric 
complex $T_L$, and how to compute some of its homotopy-type 
invariants, directly from the combinatorial data encoded in the 
input simplicial complex $L$. 

Our second---and more ambitious---goal is to understand how the 
algebraic topology of the infinite cyclic Galois covers of $T_L$ 
depends on the vertex labelings parametrizing those covers.

\subsection{Toric complexes and right-angled Artin groups} 
\label{intro:raag}

Let $L$ be a finite simplicial complex on vertex set $\V$, 
and let $T^{n}$ be the torus of dimension $n=\abs{\V}$, with 
the standard CW-decomposition. The {\em toric complex} 
associated to $L$, denoted $T_L$, is the subcomplex of 
$T^{n}$ obtained by deleting the cells corresponding to 
the non-faces of $L$. Much is known about the topology 
of these combinatorially defined spaces.  

On one hand, the fundamental group $G_L=\pi_1(T_L)$ is the 
{\em right-angled Artin group} determined by the graph 
$\G=L^{(1)}$, with presentation consisting of a generator 
$v$ for each vertex $v$ in $\V$, and a commutator relation 
$vw=wv$ for each edge $\{v,w\}$ in $\Gamma$. As shown 
by Droms \cite{Dr}, $G_L\cong G_{L'}$ if and only if the 
corresponding graphs are isomorphic.  The associated 
graded Lie algebras and the Chen Lie algebras of 
right-angled Artin groups were computed in 
\cite{DK1, DK2, PS-artin} and \cite{PS-artin}, respectively.  
For a survey of the geometric properties of such groups, 
we refer to Charney \cite{Cha}. 

On the other hand, the cohomology ring $H^*(T_L,\k)$ is the 
exterior Stanley-Reisner ring $\SR$, with generators 
the duals $v^*$, and relations the monomials corresponding 
to the missing faces of $L$, see Kim--Roush \cite{KR} 
and Charney--Davis \cite{CD}. 

\subsection{Higher homotopy groups} 
\label{intro:high pi}

As mentioned above, one of our goals here is to 
better understand the homotopy type invariants 
attached to a toric complex.  In particular, 
we seek to generalize known results, from the 
special case when $L$ is a flag complex (that is,  
one for which every subset of pairwise adjacent 
vertices spans a simplex), to the general case.  
We start in Section \ref{sect:toric} with a study 
of the higher homotopy groups of a toric complex. 

It has long been known that $T_L$ is aspherical, 
whenever $L$ is a flag complex, see \cite{CD, MV}. Recently, 
Leary and Saadeto\u{g}lu \cite{LS} have shown that the 
converse also holds. In Theorem \ref{thm:highpi}, 
we make this result more precise, by giving a 
combinatorial description of the first non-vanishing 
higher homotopy group of a non-flag complex $L$, 
viewed as a module over $\Z{G_L}$.

\subsection{Cohomology jumping loci} 
\label{intro:cjl}
The characteristic varieties $\VV^i_d(X,\k)$ and the 
resonance varieties $\RR^i_d(X,\k)$ of a finite-type 
CW-complex $X$ provide a unifying framework 
for the study of a host of questions, both quantitative and 
qualitative, concerning the space $X$ and its fundamental group.  
For instance, counting certain torsion points on the character 
torus, according to their depth with respect to the stratification
by the characteristic varieties, yields information about the 
homology of finite abelian covers of $X$. The subtle interplay 
between the geometry of these two sets of varieties leads to 
powerful formality and quasi-projectivity obstructions, 
see \cite{DPS-jump}. Finally, the cohomology jumping loci 
of a classifying space $K(G,1)$ provide computable upper 
bounds for the Bieri--Neumann--Strebel--Renz (BNSR) 
invariants of a group $G$, see \cite{PS-bns}. 

The degree $1$ jump loci of a toric complex have been 
computed in \cite{PS-artin, DPS-jump}, leading to a 
complete solution of Serre's quasi-projectivity problem 
within the class of right-angled Artin groups.  We determine 
here, in \S\ref{sect:jump toric}, the higher-degree jump loci 
of toric complexes.  Applications and further discussion 
can be found in \cite{PS-bns}. 

In Theorem \ref{thm:res toric}, we compute the 
resonance varieties $\RR^i_d(T_L,\k)\subseteq \Hom (G_L, \k)$, 
associated to the cohomology ring $\SR$, while in  
Theorem \ref{thm:cv toric}, we compute the jumping loci 
for cohomology with coefficients in rank $1$ local systems, 
$\VV^i_d(T_L,\k)\subseteq \Hom (G_L, \k^{\times})$, 
over an arbitrary field $\k$, and for all integers $i, d\ge 1$. 
Explicitly, 
\begin{equation}
\label{intro:cjl tc}
\RR^i_d(T_L,\k)= \bigcup_{\sW} \, 
\k^{\sW},
\qquad
\VV^i_d(T_L,\k)= \bigcup_{\sW} \, 
(\k^{\times})^{\sW}, 
\end{equation}
where, in both cases, the union is taken over all subsets 
$\sW\subset\sV$ for which the $i$-th ``Aomoto--Betti" 
number, $\beta_{i}(\SR,\sW)$, is at least $d$.  These  
numbers can be computed directly from $L$,  
using the following formula of Aramova, Avramov, 
and Herzog \cite{AAH}:
\begin{equation}
\label{intro:aah}
\beta_{i}(\SR,\sW)=\sum_{\sigma\in L_{\sV\setminus \sW}}
\dim_{\k} \widetilde{H}_{i-1-\abs{\sigma}} (\lk_{L_\sW}(\sigma),\k), 
\end{equation}
where $L_\sW$ is the subcomplex induced by $L$ on $\sW$, 
and $\lk_{K}(\sigma)$ is the link of a simplex $\sigma$ 
in a subcomplex $K\subseteq L$. 
In the particular case when $i=d=1$, the formulas from 
\eqref{intro:cjl tc} recover results from \cite{PS-artin} and 
\cite{DPS-jump}, respectively. 

\subsection{Infinite cyclic covers and Artin kernels}
\label{intro:art ker}
Consider a homomorphism $\chi\colon G_{L}\to \Z$,  
specified by assigning an integer weight, 
$m_v=\chi(v)$, to each vertex $v$ in $\sV$. Assume 
$\chi$ is onto, and let $\pi\colon T_L^{\chi}\to T_L$ 
be the corresponding Galois cover. The fundamental 
group 
\begin{equation}
\label{intro:ak}
N_{\chi}:=\pi_1(T_L^{\chi})=\ker (\chi\colon G_{L}\to \Z)
\end{equation}
is called the {\em Artin kernel}\/ associated to $\chi$.  
A classifying space for this group is the space 
$T_{\Delta_{\G}}^{\chi}$, where $\Delta_{\G}$ is the 
flag complex of $\G=L^{(1)}$.  As mentioned above, 
a major goal of this paper is to understand how the 
algebraic topology of the spaces $T_L^{\chi}$, and 
some of the properties of the groups $N_{\chi}$, 
depend on the epimorphism $\chi$. 

Noteworthy is the case when $\chi$ is the ``diagonal" 
homomorphism $\nu\colon G_{L}\surj \Z$, which assigns 
to each vertex the weight $1$.  The corresponding 
Artin kernel, $N_{\G}=N_{\nu}$, is called 
the {\em Bestvina--Brady group} associated to $\G$.  
As hinted at by Stallings \cite{St} and Bieri \cite{Bi} in 
their pioneering work, and as proved in full generality 
by Bestvina and Brady in their landmark paper \cite{BB}, 
the geometric and homological finiteness properties of the group 
$N_\G$ are intimately connected to the topology of the flag complex 
$\Delta_\G$. For example, $N_{\G}$ is finitely generated if and only 
if $\G$ is connected, and $N_{\G}$ is finitely presented if and only if 
$\Delta_\G$ is simply-connected.  When $\pi_1(\Delta_\G)=0$, an 
explicit finite presentation for $N_\G$ was given by Dicks and Leary 
\cite{DL}, and the presentation was further simplified in \cite{PS-bb}. 

\subsection{Homology of $\Z$-covers}
\label{intro:hom ak}

To get a handle on the coverings $T_L^{\chi} \to T_L$, we focus 
in Sections \ref{sec:hom zcovers}--\ref{sect:tm test}  
on the homology groups $H_*(T_L^{\chi},\k)$, viewed as modules 
over the group algebra $\k\Z$, with coefficients in an arbitrary 
field $\k$.   After some preparatory material in \S\ref{sec:hom zcovers}, 
we compute those homology groups, in two steps.  First, we give in 
Theorem \ref{thm:hom cov} a combinatorial formula for their 
$\k\Z$-ranks:
\begin{equation}
\label{intro:kzrank}
\rank_{\k\Z} H_{i}(T^{\chi}_L ,\k) = \beta_{i} (\SR ,\sV_0(\chi)),
\end{equation}
where $\sV_0(\chi)=\set{v \mid \chi(v)\ne 0}$ is the support of $\chi$. 

The computation of the $\k\Z$-torsion part is more complicated.  
The algorithm, which we summarize in Theorem \ref{thm:hom tors}, 
involves two steps.  The first, arithmetic  in nature, requires 
factoring certain cyclotomic polynomials in $\k[t]$. The second, 
algebro-combinatorial in nature, requires diagonalizing certain 
monomial matrices over $\k[[t]]$. In \S\ref{subsec:discuss}, 
we note that the $f$-primary part of $H_{i} (T_L^{\chi},\k)$ 
is non-trivial only when the irreducible polynomial $f\in \k [t]$ 
divides some cyclotomic polynomial, and we show that this 
is the only restriction on non-trivial torsion, for $i>0$.

As an application of our method, we give in 
Corollary \ref{cor:hom res} and Theorem \ref{thm:mono test} 
combinatorial tests for deciding whether the following 
conditions are satisfied (for fixed $r>0$):
\begin{alphenum}
\item \label{q1} For each $i\le r$, the $\k$-vector space 
$H_{i} (T_L^{\chi},\k)$ is finite-dimensional.
\item \label{q2} For each $i\le r$, the $\k\Z$-module 
$H_{i} (T_L^{\chi},\k)$ is trivial.
\end{alphenum}
The first test amounts to the vanishing of the ranks 
from \eqref{intro:kzrank}, for all $i\le r$. The second 
test amounts to the vanishing of $\beta_{i} (\SR ,\sV_{q}(\chi))$, 
for all $i\le r$, where $q$ runs through a finite list of primes 
(and $0$), depending on $\chi$ and $p=\ch\k$, and $\sV_{q}(\chi)$ 
is the support of $\chi$ in characteristic $q$.  Clearly 
\eqref{q2} $\Rightarrow$ \eqref{q1}, but not the other 
way, cf.~Remark \ref{cor:triv res}. 

For the Bestvina-Brady covers $T_L^{\nu}$, we show 
in Corollary \ref{cor:bb mono} that \eqref{q1} and \eqref{q2} 
hold simultaneously, and this happens precisely when 
$\nu_{\k} \not\in \bigcup_{i\le r} \RR_1^i(T_L, \k)$,
where $\nu_{\k}$ is the cohomology class in degree 
one naturally associated to $\nu$.

\subsection{Cohomology ring and finiteness properties}
\label{intro:coho ak}

In Section \ref{sect:coho}, we illustrate our techniques by  
determining the cohomology ring of $T_L^{\chi}$ (truncated 
in a certain degree), with coefficients in a field $\k$. 
The cohomology ring of a Bestvina-Brady cover 
$T_L^{\nu}$ was computed in \cite{LS}, up to any 
degree $r$ for which condition \eqref{q1} holds.  In 
Theorem \ref{thm:ak coho}, we extend this computation 
to arbitrary $\chi$, up to any degree $r$ for which  condition 
\eqref{q2} holds: 
\begin{equation}
\label{intro:coho}
H^{\le r}(T_L^{\chi},\k)\cong H^{\le r}(T_L,\k)/(\chik).
\end{equation}

Next, we consider the finiteness properties of Artin kernels.   
A group $G$ is said to be of type $\FP_r$ ($1\le r \le \infty$) 
if the trivial $G$-module $\Z$ has a projective $\Z{G}$-resolution, 
finitely generated in degrees up to $r$, 
see Serre \cite{Se71}. Using a result of Meier, Meinert, 
and VanWyk \cite{MMV}, as reinterpreted by Bux and 
Gonzalez \cite{BG}, we note in Theorem \ref{thm=fptriv} that 
\begin{equation}
\label{intro:fp}
\text{$N_{\chi}$ is of type $\FP_r$}
\ \Longleftrightarrow \  
\dim_{\k} H_{i} (N_{\chi}, \k)< \infty, \text{ for all $i\le r$, and 
all $\k$}.
\end{equation}
This is to be compared with a result of the same flavor 
from \cite{LS}, where it is shown that the cohomological 
dimension of a Bestvina-Brady group $N_{\nu}$ equals 
its trivial cohomological dimension.

The finiteness properties of the groups $N_{\chi}$ are 
controlled by the Bieri--Neumann--Strebel--Renz invariants 
of $G_L$, and those invariants are also computed in 
\cite{MMV, BG}. A close relationship between the first 
resonance variety $\RR^1_1(G_L,\R)$ and the BNS invariant 
$\Sigma^1(G_L)$ was first noted in \cite{PS-artin}. 
A generalization to the higher resonance varieties and the 
higher BNSR invariants, based on formula \eqref{intro:cjl tc}, 
is given in \cite{PS-bns}.

\subsection{Graded Lie algebras}
\label{intro:lie ker}

In Sections \ref{sect:holo} and \ref{sect:lie ak} we study several 
graded Lie algebras attached to groups of the form $N_{\chi}$, 
generalizing results from \cite{PS-bb}, valid only in the case 
where $\chi=\nu$, the diagonal character. 

We start with some general properties of holonomy 
Lie algebras. Given a graded algebra $A$ satisfying  
some mild assumptions, its {\em holonomy Lie algebra}, 
$\h(A)=\bigoplus_{s\ge 1} \h_s(A)$, is defined as 
the free Lie algebra on the dual of $A^1$, modulo 
the (homogeneous) ideal generated by the image of the 
comultiplication map. The main result here is 
Theorem \ref{thm:hab}, where we show that, if an 
element $a\in A^1$ is non-resonant up to degree $r$, then 
\begin{equation}
\label{intro:hab}
\left(\h(A)/\h''(A)\right)_{s} \cong\left(\h(A/aA)/\h''(A/aA)\right)_{s}, 
\quad \text{for $2\le s\le r+1$}.
\end{equation}

Next, we consider the associated graded Lie algebra, 
$\gr(N_{\chi})$, arising from the lower central series of 
an Artin kernel $N_{\chi}$, and the rational holonomy 
Lie algebra, $\h(N_{\chi})$, arising from the cohomology 
ring $A=H^*(N_{\chi}, \Q)$. In Proposition \ref{prop:gr split} 
and Corollary \ref{cor:gr ak}, we assume that $H_1(N_{\chi}, \Q)$ 
is $\Q\Z$-trivial, and determine the associated graded Lie 
algebra, $\gr(N_{\chi})$, and the Chen Lie algebra, 
$\gr(N_{\chi}/N''_{\chi})$. Using computations from 
\cite{PS-artin}, we show that their graded ranks, 
$\phi_k$ and $\theta_k$, are given by 
\begin{equation}
\label{intro:clique cut}
\prod_{k=1}^{\infty}(1-t^k)^{\phi_k}=\frac{P(-t)}{1-t} 
\quad\text{and}\quad 
\sum_{k=2}^{\infty} \theta_k t^{k} = Q \Big(\frac{t}{1-t}\Big),
\end{equation}
where $P(t)$ and $Q(t)$ are the clique and cut polynomials 
of the graph $\G$. In Theorem \ref{thm:holo split}, we 
assume additionally that $H_2(N_{\chi}, \Q)$ is $\Q\Z$-trivial, 
and determine the graded Lie algebra $\h(N_{\chi})$ in this 
instance; in particular, we find that $\h'(N_{\chi})=\h'(G_{\G})$. 

\subsection{The $1$-formality property}
\label{intro:formal}

A finitely generated group $G$ is said to be {\em $1$-formal}, 
in the sense of Sullivan \cite{S}, if there is a filtration-preserving 
Lie algebra isomorphism between the Malcev Lie algebra 
of $G$, as defined by Quillen \cite{Q}, and the degree 
completion of the rational holonomy Lie algebra of $G$, 
as defined by Chen \cite{Ch77}. Put another way, 
$G$ is $1$-formal if its Malcev Lie algebra, $\M(G)$, 
is quadratically presented.  Examples include fundamental 
groups of compact K\"{a}hler manifolds and complements 
of algebraic hypersurfaces in $\CP^n$, certain pure braid groups of 
Riemann surfaces, and certain Torelli groups.  The $1$-formality 
property of a group has many remarkable consequences.  
We refer to \cite{S, PS-chen, DPS-jump} for more details 
and references. 

In \cite{KM}, Kapovich and Millson showed that all finitely 
generated Artin groups (in particular, all groups of the 
form $G_\G$) are $1$-formal. In \cite{PS-bb}, we showed 
that all finitely presented Bestvina-Brady groups, i.e., all 
groups $N_{\G}$ for which $\pi_1(\Delta_\G)=0$, are $1$-formal. 

We generalize this last result here, in Theorem \ref{thm:ker formal}: 
Suppose $H_i(N_{\chi}, \Q)$ is $\Q\Z$-trivial, for $i=1,2$
(when $\chi=\nu$, this  
boils down to $\widetilde{H}_{i}(\Delta_{\G}, \Q)=0$, for 
$i=0,1$); then the Artin kernel $N_{\chi}$ is a finitely generated, 
$1$-formal group. Using this Theorem, we are able to construct 
$1$-formal groups which are not finitely presentable.  To the 
best of our knowledge, these examples are the first of their kind.

\section{Toric complexes and right-angled Artin groups}
\label{sect:toric}

In this section, we discuss in more detail some of our 
main characters---the spaces $T_L$ and their fundamental 
groups, $G_L$---and compute certain homotopy groups of $T_L$. 

\subsection{Toric complexes}
\label{subsec:tl}
Let $L$ be a finite simplicial complex, on vertex set $\sV$. 
For each simplex $\sigma=\set{v_{i_1},\dots,v_{i_k}}$ of $L$, 
let $T_{\sigma}$ be the torus formed by identifying parallel 
faces of a $k$-cube ($T_{\emptyset}$ is a point).  
The {\em toric complex} associated to $L$ is the 
identification space 
\begin{equation}
\label{eq:tc}
T_L=\coprod_{\sigma \in L} T_{\sigma} 
\Big\slash \big\{ 
T_{\sigma} \cap T_{\sigma'} = T_{\sigma\cap\sigma'}
\big\}\Big. .
\end{equation}

An equivalent description is as follows.  Let $T^{n}$ be the 
torus of dimension $n=\abs{\V}$, with the standard 
CW-decompo\-sition.   Then $T_L$ is the subcomplex of 
$T^{n}$ obtained by deleting the cells corresponding to 
the non-faces of $L$. Note that the $k$-cells $c_\sigma$ 
in $T_L$ are in one-to-one correspondence with the 
$(k-1)$-simplices $\sigma$ in $L$. 

In the terminology from \cite{DS07}, $T_L=\mathcal{Z}_L(S^1)$ 
is an example of a ``generalized moment-angle complex."  As such, 
the construction enjoys many natural properties.  For instance, 
$T_{K*L}=T_{K}\times T_{L}$, where $*$ denotes simplicial join; 
moreover, if $K\subset L$ is a simplicial subcomplex, then 
$T_K\subset T_L$ is a CW-subcomplex.

\subsection{The cohomology ring}
\label{subsec:sr}
Fix a coefficient ring $\k$ (we will mainly be interested in the 
case when $\k=\Z$, or $\k$ is a field). Let $C_{\bullet}(L,\k)$ 
be the simplicial chain complex of $L$, and let $C_{\bullet}(T_L,\k)$ 
be the cellular chain complex of $T_L$. Since $T_L$ is a 
subcomplex of $T^{n}$, all differentials in $C_{\bullet}(T_L,\k)$ 
vanish, i.e., $T_L$ is a {\em minimal} CW-complex. It follows that 
\begin{equation}
\label{eq:hom tl}
H_k(T_L,\k)=C_{k-1}(L,\k),\quad \text{for all $k>0$}. 
\end{equation}
In other words, $H_k(T_L,\k)$ is a free $\k$-module of 
rank $d_k(L)=\#\set{\sigma \in L\mid \abs{\sigma}=k}$, 
where $\abs{\sigma}=\dim(\sigma)+1$.  In particular, the 
Betti numbers $b_k(T_L)=\dim_{\k} H_k(T_L,\k)$ depend 
only on $L$, and not on $\k$. 

As shown in \cite{KR, CD}, the cohomology ring of $T_L$ 
may be identified with the exterior Stanley--Reisner ring of $L$.  
More precisely, let $V$ be the free $\k$-module on the set 
$\V$, and $V^*=\Hom_{\k}(V,\k)$ its dual.  Set 
$\SR=\bigwedge V^* /J_L$, 
where $\bigwedge V^*$ is the exterior algebra on $V^*$ 
and $J_L$ is the ideal  generated by all monomials 
$t_{\sigma}=v^*_{i_1} \cdots v^*_{i_k}$ corresponding 
to simplices $\sigma=\set{v_{i_1}, \ldots,v_{i_k}}$ not belonging 
to $L$. Then 
\begin{equation}
\label{eq:sr ring}
H^*(T_L,\k)=\SR.
\end{equation}
In the case when $L$ is a flag complex, $J_L$ is a quadratic 
monomial ideal, and so $\SR$ is a Koszul algebra. 

\subsection{Right-angled Artin groups}
\label{subsec:raag}
The $1$-skeleton $L^{(1)}$ can be viewed as a (finite) simple 
graph $\G=(\V,\E)$, with vertex set $\V$ consisting of the 
$0$-cells, and edge set $\E$ consisting of the $1$-cells. 
It is readily seen that the fundamental group $G_L=\pi_1(T_L)$ 
is isomorphic to 
\begin{equation}
\label{eq:artin group}
G_{\G}:= \langle v \in \V \mid 
uv=vu \text{ if $\{u,v\} \in \E$}\rangle, 
\end{equation}
the right-angled Artin group associated to $\G$. 
These groups interpolate between free groups 
(if $\G=\overline{K}_n$ is the discrete graph on 
$n$ vertices, then $G_{\G}=F_n$), and free abelian 
groups (if $\G=K_n$ is the complete graph on $n$ 
vertices, then $G_{\G}=\Z^n$). Moreover, if 
$\G=\G' \coprod \G''$ is the disjoint union of two graphs, 
then $G_{\G}= G_{\G'}*G_{\G''}$, and if $\G=\G' * \G''$ 
is their join, then $G_{\G}= G_{\G'}\times G_{\G''}$.

Let  $\Delta=\Delta_\G$ be the flag complex of $\G$, i.e., 
the maximal simplicial complex with $1$-skeleton equal to $\G$.  
Clearly, $L$ is a subcomplex of $\Delta$, sharing the same 
$1$-skeleton.  Moreover, the $k$-simplices of $\Delta$ correspond 
to the $(k+1)$-cliques of $\G$.  A classifying space for $G=G_\G$ 
is the toric complex associated to $\Delta$,
\begin{equation}
\label{eq:kg1}
K(G,1)=T_{\Delta}, 
\end{equation}
see Charney--Davis \cite{CD} and Meier--VanWyk \cite{MV}. 
For example, if $\G=\overline{K}_n$, then $\Delta=\G$, 
and $T_{\Delta}=\bigvee^n S^1$, whereas if $\G=K_n$, 
then $\Delta=\Delta_{n-1}$, the $(n-1)$-simplex, and $T_{\Delta}=T^n$. 

Let $\widetilde{C}_{\bullet}=(C_{\bullet}(\widetilde{T}_{L},\k),
\tilde\partial_{\bullet})$ be the equivariant chain complex of the 
universal cover of $T_L$.  Under the identification 
$\widetilde{C}_k = \k G \otimes_{\k} C_k$, where 
$C_k=C_k(T_L,\k)$, 
the boundary map $\tilde\partial_k\colon\, 
\widetilde{C}_k\to \widetilde{C}_{k-1}$ is given by
\begin{equation}
\label{eq:equiv bdry}
\tilde\partial_k (1\otimes c_{\sigma}) = \sum_{r=1}^{k}  
(-1)^{r-1} (v_{r} -1) \otimes c_{\sigma\setminus \{v_r\} }
\end{equation}
where $\sigma=\{v_{1},\dots , v_{k}\}$ is a $(k-1)$-simplex 
in $L$. 

\subsection{Higher homotopy groups of toric complexes}
\label{subsect:highpi}

By the above-mentioned result, all the higher homotopy 
groups of the toric complex associated to a flag complex 
vanish.  In general, though, $T_L$ will 
have $\pi_r(T_L) \ne 0$, for some $r>1$.   
The next result identifies the first integer $r$ 
for which this happens, and computes the 
rank of the $G_L$-coinvariants for the corresponding 
homotopy group (viewed as a module over $\Z{G_L}$). 

\begin{theorem}
\label{thm:highpi}
For a finite simplicial complex $L$, with  
associated flag complex $\Delta$, set 
$p=p(L):=\sup\, \{k \mid d_{\le k}(\Delta)=d_{\le k}(L)\}$. 
Let $T_L$ be the corresponding toric complex, 
and $G=\pi_1(T_L)$. Then, the following 
are equivalent:
\begin{romenum}
\item \label{pi1} $T_L$ is aspherical. 
\item \label{pi2} $L$ is a flag complex.
\item \label{pi3} $p=\infty$. 
\end{romenum}
Moreover, if $p<\infty$, then:
\begin{romenum}
\setcounter{enumi}{3}
\item \label{pi4}
$\pi_2(T_L)=\cdots =\pi_{p-1}(T_L)=0$. 
\item \label{pi5}
The $p$-th homotopy group of $T_L$, when viewed 
as a module over $\Z{G}$, has a finite presentation  
of the form $\pi_{p}(T_L) = 
\coker (\Pi_{p+1} \circ \tilde\partial_{p+2})$, 
where $\Pi_{p+1} \colon 
\Z{G}^{d_{p+1}(\Delta)} \to 
\Z{G}^{d_{p+1}(\Delta)-d_{p+1}(L)}$ 
is the canonical projection.  
\item \label{pi6} 
The group of coinvariants $\big(\pi_{p}(T_L)\big)_G$  
is free abelian, of rank $d_{p+1}(\Delta)-d_{p+1}(L)>0$. 
\end{romenum}
\end{theorem}

\begin{proof}
First suppose $p<\infty$. Then $L$ is a proper simplicial 
subcomplex of $\Delta$, and both have the same $1$-skeleton, 
$\Gamma$. Thus, $T_L$ is a proper subcomplex of 
the aspherical CW-complex $T_\Delta$, and both share 
the same $2$-skeleton, $T_{\Gamma}$. Moreover, 
$T_L$ and $T_\Delta$ are minimal CW-complexes, 
and their cohomology rings, $H^*(T_L;\Z)= 
\Z\langle T_L\rangle$ and $H^*(T_\Delta;\Z)= 
\Z\langle T_\Delta\rangle$, are generated in degree $1$. 
Statements \eqref{pi4}--\eqref{pi6} now follow at once 
from \cite{PS-highpi}, Theorem 2.10, Corollary 2.11, 
and Remark 2.13.  

The implication \eqref{pi3} $\Rightarrow$ \eqref{pi2} is 
clear, whereas \eqref{pi2} $\Rightarrow$ \eqref{pi1} follows 
from \eqref{eq:kg1}.  To prove  \eqref{pi1} $\Rightarrow$ \eqref{pi3}, 
assume $p$ is finite. Then $\pi_{p}(T_L)\ne 0$, by \eqref{pi6};  
thus, $T_L$ is not aspherical. 
\end{proof}

Equivalence \eqref{pi1} $\same$ \eqref{pi2} in the above recovers 
a result of Leary and Saadeto\u{g}lu (\cite[Proposition 4]{LS}), 
which they proved by different means. 

\section{Aomoto-Betti numbers and cohomology jumping loci}
\label{sect:jump toric}

In this section, we determine the resonance varieties 
$\RR^i_d(T_L, \k)$ and characteristic varieties 
$\VV^i_d(T_L, \k)$  of a toric complex $T_L$, in all 
degrees $i,d\ge 1$, and for all fields $\k$. 

\subsection{Resonance varieties of algebras}
\label{subsec:res var}

Let $A$ be a connected, locally finite, graded algebra over 
a field $\k$.  Denote the Betti numbers of $A$ by 
$b_{i}(A)=\dim_{\k} A^i$.  For each element $z\in A^1$ 
with $z^2=0$, right-multiplication by $z$ defines a 
cochain complex 
\begin{equation}
\label{eq:aomoto}
\xymatrixcolsep{20pt}
\xymatrix{(A , z)\colon \
A^0 \ar^(.65){z}[r] & A^1 \ar^(.48){z}[r] 
& \cdots  \ar[r] & A^{i-1}  \ar^(.53){z}[r] & A^i \ar^(.42){z}[r] 
& A^{i+1}  \ar[r] & \cdots},
\end{equation}
also known as the {\em Aomoto complex}. 
Denote the Betti numbers of this complex by 
\begin{equation}
\label{eq:aomoto betti}
\beta_{i}(A,z)=\dim_{\k} H^i(A, z).
\end{equation} 
\begin{remark}
\label{rem:beta}
The following properties of the Aomoto Betti numbers 
are immediate:
\begin{enumerate}
\item $\beta_{i}(A,z)\le b_{i}(A)$, for all $z\in A^1$, 
and $\beta_{i}(A,0)=b_{i}(A)$.
\item $\beta_0(A,z)=0$, for all $z\ne 0$, 
and $\beta_0(A,0)=1$. 
\end{enumerate}
\end{remark}

Suppose $z^2=0$, for all $z\in A^1$.  Then, for each 
$i\ge 1$ and each $0\le d \le b_{i}=b_{i}(A)$, we may define 
the {\em resonance variety}
\begin{equation}
\label{eq:res var alg}
\RR^i_d(A)=\{z \in A^1 \mid \beta_{i}(A,z) \ge  d\}.
\end{equation}

It is readily seen that each set $\RR^i_d=\RR^i_d(A)$ is a 
homogeneous algebraic subvariety of the affine 
space $A^1=\k^{b_1}$. For each degree $i\ge 1$, the 
resonance varieties provide a descending filtration 
$A^1=\RR^i_0 \supset \RR^i_1 \supset \cdots 
\supset \RR^i_{b_{i}} \supset \RR^i_{b_{i+1}}=\emptyset$.

\begin{remark}
\label{rem:beta vanish}
Given $z\in A^1$, and an integer 
$r\ge 1$, the following are equivalent:
\begin{enumerate} 
\item \label{br1}  $\beta_{i} (A,z)=0$, for all $i\le r$.
\item \label{br2}  $z\notin \RR^i_1(A)$, for all $i\le r$. 
\end{enumerate}
\end{remark}

A simple example is provided by the exterior algebra 
$E=\bigwedge \k^n$.  It is readily checked that the Aomoto 
complex \eqref{eq:aomoto} is exact in this case.  Hence, 
$\RR^i_d(E)=\{0\}$, for all $1\le i\le n$ and $1\le d\le \binom{n}{i}$.  

Now let $X$ be a connected CW-complex with finitely many 
cells in each dimension.  Fix a field $\k$, and let $A=H^* (X,\k)$ 
be the respective cohomology ring.  (If $\ch\k=2$, we need to 
assume $H_1(X,\Z)$ is torsion-free.) Then $A$ is a locally finite, 
graded algebra with $z^2=0$, for all $z\in A^1$.  
Thus, for every $z\in A^1$, we may define the cochain 
complex $(A,\cdot z)$, with Aomoto Betti numbers 
$\beta^{\k}_{i}(X, z) := \beta_{i} (H^*( X,\k) , z)$. 
The resonance varieties of $X$ (over the field $\k$) 
are simply those of its cohomology ring:
$\RR^i_d(X,\k):=\RR^i_d(H^*(X,\k))$.  Similarly, 
if $G$ is a group admitting a classifying space of
finite type, $\RR^i_d(G,\k):=\RR^i_d(H^*(G,\k))$.
 
\subsection{Aomoto Betti numbers for $\SR$}
\label{subsec:aomoto sr}

We now specialize to the case when $X=T_L$ is the 
toric complex associated to a finite simplicial complex 
$L$, and $\SR=H^*(T_L,\k)$ is the corresponding exterior 
face ring. Since $H_1(T_L,\Z)$ is torsion-free, the resonance 
varieties $\RR^i_d(T_L,\k)=\RR^i_d(\SR)$ are defined over 
any field $\k$.

From the definitions, $\SR^1=V^*$, where $V$ is the $\k$-vector 
space with basis indexed by the vertex set $\sV$ of $L$, and 
$V^*$ its dual.   Given an element $z\in V^*$, write 
$z=\sum_{v\in \sV} z_v v^*$, with $z_v\in \k$, and 
define the {\em support}\/ of $z$ as 
\begin{equation}
\label{eq:supp}
\supp (z)=\set{v \in \V\mid z_v \ne 0}.  
\end{equation}

Conversely, given a subset $\sW \subset \sV$ of the vertex set, 
we may define a ``canonical" element $z_{\sW}\in V^*$ by 
\begin{equation}
\label{eq:zw}
z_{\sW}=\sum_{v\in \sW} v^*, 
\end{equation}
with the convention that $z_{\emptyset}=0$. 
Obviously, $\supp(z_{\sW})=\sW$.  For simplicity, we 
will write 
\begin{equation}
\label{eq:beta w}
\beta_{i}(\SR,\sW):=\beta_{i}(\SR,z_{\sW}).
\end{equation} 

\begin{remark}
\label{rem:ab sr}
Unlike the Betti numbers $b_{i}(T_L)=\dim_{\k} H_{i}(T_L,\k)$, 
the Aomoto Betti numbers $\beta_{i}(\SR,\sW)$ do depend 
on $\k$ (in fact, just on $p=\ch \k$), and not only on $L$ (and 
$\sW$); see Proposition \ref{prop:hochster aah} below. 
Note also that $\beta_0(\SR,\sW)=0$, whenever 
$\sW\ne \emptyset$, and  $\beta_{i}(\SR,\emptyset)=b_{i}(T_L)$. 
\end{remark} 

The following result is due to Aramova, Avramov, 
and Herzog (Proposition 4.3 in \cite{AAH}). For 
the sake of completeness, we reproduce the proof. 

\begin{lemma}[\cite{AAH}] 
\label{lem:aah}
Let $z$, $z'\in \SR^1$. If  $\supp(z)=\supp(z')$,  
then $\beta_{i}( \SR , z) = \beta_{i}( \SR , z')$, 
for all $i\ge 1$. 
\end{lemma}

\begin{proof}
Write $z=\sum_{v\in \sW} z_v v^*$ and 
$z'=\sum_{v\in \sW} z'_v v^*$, where $\sW$ 
is the common support.  The linear map 
$\phi\colon V^*\to V^*$ given on basis 
elements by 
\begin{equation}
\label{eq:phi map}
\phi(v^*)=\begin{cases}
\frac{z'_v}{z_v} v^* &\text{for $v\in \sW$}\\ 
v^* &\text{otherwise},
\end{cases}
\end{equation}
extends to an algebra isomorphism 
$\phi\colon \SR\to \SR$, taking $z$ to $z'$. 
The conclusion follows. 
\end{proof}

\begin{lemma}
\label{lem:semi}
If $\sW'\subset \sW$, then 
$\beta_{i}(\SR,\sW')\ge \beta_{i}(\SR,\sW)$, for all $i\ge 1$.
\end{lemma}
\begin{proof}
For each $t\in \k$, define 
\[
z_t=\sum_{v\in \sW'} v^* + t \sum_{v\in \sW\setminus \sW'} v^*. 
\]
Note that $z_0=z_{\sW'}$ and $z_1=z_{\sW}$. By lower  
semicontinuity of Betti numbers,
\[
\beta_{i}(\SR, z_t) \le \beta_{i}(\SR,z_0),
\]
for $t$ in a Zariski open subset of $\k$, containing $0$. 
On the other hand, for $t\ne 0$, $\beta_{i}(\SR, z_t)=\beta_{i}(\SR, z_1)$, 
by Lemma \ref{lem:aah}. This finishes the proof. 
\end{proof}

\subsection{The Aramova--Avramov--Herzog formula}
\label{subsec:aah formula}
The Aomoto Betti numbers of the exterior face ring $\SR$ 
can be computed in purely combinatorial terms, by using 
the following Hochster-type formula from \cite[Proposition 4.3]{AAH}, 
suitably interpreted and corrected. 

\begin{prop}[\cite{AAH}] 
\label{prop:hochster aah}
Let $L$ be a finite simplicial complex on vertex set $\sV$, 
and let $\sW\subset \sV$. Then:
\[
\beta_{i}(\SR,\sW)=\sum_{\sigma\in L_{\sV\setminus \sW}}
\dim_{\k} \widetilde{H}_{i-1-\abs{\sigma}} (\lk_{L_\sW}(\sigma),\k).
\]
\end{prop}
Here $L_\sW=\set{\tau\in L\mid \tau \subset \sW}$ is the 
simplicial complex obtained by restricting $L$ to $\sW$, and 
$\lk_{L_\sW}(\sigma)= \set{\tau \in L_\sW \mid \tau\cup \sigma \in L}$ 
is the link of a simplex 
$\sigma$ in $L_\sW$. The range of summation 
in the above formula includes the empty simplex, 
with the convention that $\abs{\emptyset}=0$ and 
$\widetilde{H}_{-1}(\emptyset,\k)=\k$. 
In particular, 
\begin{align}
\label{eq:beta1}
\beta_1(\SR,\sW)
&=\dim_{\k} 
\widetilde{H}_{0} (\lk_{L_\sW}(\emptyset),\k)
+\sum_{v\in \sV \setminus \sW} \dim_{\k} 
\widetilde{H}_{-1} (\lk_{L_\sW}(v),\k)\\
&=\tilde{b}_0(L_\sW)+ 
\abs{\set{v\in \sV \setminus \sW \mid \lk_{L_\sW}(v)= \emptyset}}. \notag
\end{align}

From the proposition, we obtain the following immediate corollary.

\begin{corollary}
\label{cor:beta vanish}
Let $L$ be a finite simplicial complex on vertex set $\sV$.   
For a subset $\sW\subset \sV$, and an integer $r>0$, 
the following are equivalent:
\begin{romenum}
\item \label{ah1} $\beta_{i}(\SR,\sW)=0$, for $1\le i\le r$. 
\item \label{ah2}  
$\widetilde{H}_{i} (\lk_{L_\sW}(\sigma),\k)=0$, 
for  all $\sigma\in L_{\sV\setminus \sW}$ and 
$-1\le i\le r-1-\abs{\sigma}$. 
\end{romenum}
\end{corollary}

\subsection{Resonance of toric complexes}
\label{subsec:res tk}
Denote the $\k$-vector space $V^*=\SR^{1}$ by 
$\k^{\sV}$.  For a subset $\sW \subset \sV$, denote 
by $\k^\sW$ the corresponding coordinate subspace. 

\begin{theorem}
\label{thm:res toric}
Let $L$ be a finite simplicial complex on vertex set $\sV$. 
Then, the resonance varieties of the toric complex $T_L$ 
(over a field $\k$) are given by:
\[
\RR^i_d(T_L,\k)= \bigcup_{\stackrel{\sW \subset \sV}%
{\beta_{i}(\SR,\sW)\ge d}} \k^ \sW.
\]
\end{theorem}

\begin{proof}
Suppose $z\in \RR^i_d(T_L,\k)$.  
By definition, this means $\beta_{i}( \SR , z) \ge d$.   
Set $\sW=\supp (z)$; then $z$ belongs to the subspace $
\k^{\sW}$.  By Lemma \ref{lem:aah}, $\beta_{i}( \SR, \sW) \ge d$. 
 
Now suppose $z\in \k^{\sW}$, 
for some subset $\sW \subset \sV$ for which  
$\beta_{i}(\SR,\sW)\ge d$. Write $\sW'=\supp (z)$; 
clearly, $\sW'\subset \sW$.  By Lemmas \ref{lem:aah} 
and \ref{lem:semi}, 
\[
\beta_{i}(\SR,z)= \beta_{i}(\SR,\sW')\ge \beta_{i}(\SR,\sW)\ge d,
\]
and so $z\in \RR^i_d(T_L,\k)$.
\end{proof}

As a corollary, we recover the description from 
\cite[Theorem 5.5]{PS-artin} of the first resonance variety 
of a right-angled Artin group. 

\begin{corollary}[\cite{PS-artin}]
\label{thm:res artin}
Let $\G=(\V,\E)$ be a finite graph.  Then 
$\RR^1_1(G_{\G},\k) = \bigcup_{\sW} \k^{\sW}$,  
where the union is over all subsets $\sW \subset \sV$ 
such that the induced subgraph $\Gamma_{\sW}$ is 
disconnected. 
\end{corollary}

\begin{proof}
By definition, the resonance varieties of $G_{\G}$ are those 
of $\SR$, where $L=\Delta_{\Gamma}$. By \eqref{eq:beta1}, 
we have $\beta_1(\SR,\sW)=0$ if and only if $L_{\sW}$ is 
connected and dominating, i.e., for all $v\in \sV\setminus \sW$, 
there is a $w\in \sW$ such that $\set{v,w}\in \sE$. The conclusion 
easily follows. 
\end{proof}

\subsection{Non-propagation of resonance}
\label{subsec:res nonprop}
One may wonder whether, for a graded algebra $A$ as in 
\S\ref{subsec:res var}, resonance ``propagates," i.e., whether 
$z\in \RR^{i}_1(A)$ implies $z\in \RR^{k}_1(A)$, for all $k\ge i$ 
such that $A^j \ne 0$, for $i\le j\le k$.  Such a phenomenon 
is believed to hold when $A$ is the Orlik-Solomon algebra of 
a complex hyperplane arrangement. The following example 
shows that resonance in degree $1$ does not propagate in 
higher degrees, even for the exterior Stanley-Reisner rings 
of flag complexes.

\begin{example}
\label{ex:nonprop}
Let $\G=\G_1 \coprod \G_2$, where $\G_{j}=K_{n_{j}}$ 
are complete graphs on $n_j\ge 2$ vertices ($j=1,2$),  
and consider the toric complex $T_{\Delta_\G}$.
The simplest example (with $n_1=n_2=2$) is the 
graph $\G=\Delta_\G$ depicted below:
\begin{equation*}
\label{eq:two edges}
\xygraph{
[]!{(0,0)}
*{\cir<2pt>{}}
-@[]!{(1,0)}
*{\cir<2pt>{}}
[]!{(2,0)}
*{\cir<2pt>{}}
-@[]!{(3,0)}
*{\cir<2pt>{}}
}
\end{equation*}
Using Theorem \ref{thm:res toric}, it is readily seen that
\[
\RR_1^i (T_{\Delta_\G}, \k)= 
\begin{cases}
\k^{n_1+n_2}, &\text{if $i=1$},\\
\k^{n_1}\times \set{0} \cup \set{0}\times \k^{n_2}, 
&\text{if $1<i\le \min (n_1, n_2)$} .
\end{cases}
\]
\end{example} 

\subsection{Characteristic varieties}
\label{subsec:cv toric}

Let $X$ be a connected CW-complex with finitely many 
cells in each dimension. Let $G=\pi_1(X)$ be the fundamental 
group, and $\Hom(G,\k^{\times})$ its group of $\k$-valued 
characters.  The {\em characteristic varieties} of $X$ 
(over $\k$) are the jumping loci for homology with coefficients 
in rank~$1$ local systems:
\begin{equation}
\label{eq:cvs}
\VV^i_d(X,\k)=\{\rho \in \Hom(G,\k^{\times}) 
\mid \dim_\k H_{i}(X,\k_{\rho})\ge d\}. 
\end{equation}
Here $\k_{\rho}$ is the $1$-dimensional $\k$-vector 
space, viewed as a module over the group ring $\k{G}$ 
via $g \cdot a = \rho(g)a$, for $g\in G$ and $a\in \k$. 
By definition, $H_*(X,\k_{\rho})$ is the homology of 
the chain complex 
$(\k_{\rho} \otimes_{\k{G}} C_{\bullet}(\widetilde{X}, \k),
\partial^{\rho})$, where 
$\widetilde{X}$ is the universal cover, 
$(C_{\bullet}(\widetilde{X}, \k),\tilde{\partial})$ is the 
($G$-equivariant) cellular chain complex, and 
$\partial^{\rho}=\id \otimes_{\k{G}} \tilde\partial$. 

Now let $X=T_L$ be the toric complex corresponding 
to a finite simplicial complex $L$ on vertex set $\sV$. 
The character variety $\Hom(G_L,\k^{\times})$ may be identified 
with the algebraic torus $(\k^{\times})^{\sV}$.  
For a subset $\sW \subset \sV$, denote by 
$(\k^{\times})^{\sW}$ the corresponding 
coordinate subtorus. 

\begin{lemma}
\label{lem:change vars}
Let $\rho\colon G_{L}\to \k^{\times}$ be a character.  Then 
\[
\dim_{\k} H_{i}(T_L,\k_{\rho}) = \beta_{i}(\SR, z_{\rho}),
\]
where $z_{\rho}=\sum_{v\in \sV}(\rho(v)-1) v^*$. 
\end{lemma}

\begin{proof}
The equivariant chain complex 
$(C_{\bullet}(\widetilde{T}_{L},\k),\tilde\partial_{\bullet})$ 
has boundary maps given by \eqref{eq:equiv bdry}. 
Note that 
$\partial^{\rho}_k (1\otimes c_{\sigma}) = \sum_{r=1}^{k}  
(-1)^{r-1} (\rho(v_{r}) -1) \otimes c_{\sigma\setminus \{v_r\} }$. 
It is readily checked that the chain complex 
$(\k_{\rho} \otimes_{\k{G}} C_{\bullet}(\widetilde{X}, \k),
\partial^{\rho})$ is dual to the cochain complex 
$(\SR,z_{\rho})$.  
\end{proof}

Using Theorem \ref{thm:res toric} and Lemma 
\ref{lem:change vars}, an argument similar to 
that in  \cite[Proposition 10.5]{DPS-serre} yields 
the following description of the characteristic 
varieties of a toric complex. 

\begin{theorem}
\label{thm:cv toric}
Let $L$ be a finite simplicial complex on vertex set $\sV$. 
Then, 
\[
\VV^i_d(T_L,\k)= \bigcup_{\stackrel{\sW \subset \sV}%
{\beta_{i}(\SR,\sW)\ge d}}
(\k^{\times})^{\sW}.
\]
\end{theorem}

As a corollary, we recover the description from 
\cite[Proposition 10.5]{DPS-serre} of the first 
characteristic variety of a right-angled Artin group. 

\begin{corollary}
\label{cor:cv raag}
Let $\G=(\V,\E)$ be a finite graph.  Then 
$\VV^1_1(G_{\G},\k) = \bigcup_{\sW} (\k^{\times})^{\sW}$,  
where the union is over all subsets $\sW \subset \sV$ 
such that the induced subgraph $\Gamma_{\sW}$ is 
disconnected. 
\end{corollary}

\section{Homology of infinite cyclic covers}
\label{sec:hom zcovers}

Let $X^{\nu}\to X$ be a Galois $\Z$-cover. In this section, 
we study the homology groups $H_*(X^{\nu},\k)$, 
viewed as modules over $\k\Z$, and the cohomology ring 
$H^*(X^{\nu},\k)$.  We start with some of the relevant 
algebraic background. 

\subsection{Finitely generated modules over $\k\Z$}
\label{subsec:laurent}
Let $\k$ be a field.  We will identify the group algebra 
$\k\Z$ with the ring of Laurent polynomials, $\L=\k[t^{\pm 1}]$.  
The irreducible polynomials in $\k[t^{\pm 1}]$ coincide (up 
to units) with the irreducible polynomials, different from $t$, 
in the subring $\k[t]$. Since $\L$ is a principal ideal domain, 
every finitely generated $\L$-module $M$ decomposes 
as a finite direct sum 
\begin{equation}
\label{eq:pid}
M=F(M) \oplus \bigoplus_{\stackrel{t\ne f\in \k[t]}%
{f\:\text{irreducible}}}  T_f(M),
\end{equation}
where $F(M)=\Lambda^{\rank M}$ denotes the free part, 
and $T_f(M)=\bigoplus_{j\ge 1} (\L/f^j \L)^{e_{j,f}(M)}$
denotes the $f$-primary part. 

Particularly simple is the case when $f=t-1$ 
and $j=1$. Note that $\L/(t-1)\L=\k$, with module structure given 
by $t\cdot 1=1$.  We say a  finitely generated $\L$-module $M$ 
is {\em trivial} (or, has trivial $\Z$-action), if $M$ decomposes 
as a direct sum of such modules, i.e., $M=(\L/(t-1)\L)^e$; 
equivalently, $t \cdot m=m$, for all $m\in M$. 

Suppose now $\k$ is algebraically closed.  Then all irreducible 
polynomials in $\L$ are of the form $t-a$, for some 
$a\in \k^{\times}$. For simplicity, write 
\begin{equation}
\label{eq:tam}
T_a(M):=T_{t-a}(M)=\bigoplus_{j\ge 1} 
(\Lambda/(t-a)^j \Lambda)^{e_{j,a}(M)}.
\end{equation}
The homomorphism $\ev_a\colon \Z\to \k^{\times}$, 
$\ev_a(n)= a^n$ defines a $\Lambda$-module $\k_a$, 
via $g \cdot x:=g(a)x$, for $g\in \Lambda$ and $x\in \k$.  
Clearly, $\k_a=T_a(\k_a)=\Lambda/(t-a) \Lambda$. 

\subsection{The $(f)$-adic completion of $\k\Z$}
\label{subsec:completion}
Returning to the general case, 
fix an irreducible polynomial $f\in \k[t]$, $f\ne t$, and let 
$\wL$ be the completion of $\L=\k[t^{\pm 1}]$ with respect 
to the filtration $(f^k)_{k\ge 0}$. It is easy to see that 
the $(f)$-adic filtration on $\L$ is separated. Hence, 
the canonical map $\iota\colon \L\to \wL$ is injective. 

Denote by $(\widehat{f^k})_{k\ge 0}$ the corresponding 
filtration on $\wL$.  By construction, $\iota\colon \L\to \wL$ 
preserves filtrations, and induces an isomorphism 
$\gr(\iota)\colon \gr(\L)\to \gr(\wL)$ between 
the associated graded rings. Note that $\gr^0(\L)$ 
is simply the residue field $\k_f:=\k[t]/(f)$. 

\begin{remark}
\label{rem:units}
Let $u\in \k[t]$ be a polynomial, $\iota(u)$ its image in $\wL$, 
and $\bar{u}$ its image in $\k_f$. Then $\iota(u)$ is a unit in 
$\wL$ $\same$ $\bar{u}\ne 0$ in $\k_f$ $\same$ $f \nmid u$. 
\end{remark}

Denote by $\lambda_f$ the homothety associated to $f$.

\begin{lemma}
\label{lem:gr iso}
There is an isomorphism of graded $\k_f$-algebras, 
$\gr(\wL)\cong \k_f[t]$.
\end{lemma}

\begin{proof}
We may replace $\gr(\wL)$ by $\gr(\L)$. It is straightforward to 
check that $\lambda_f$ induces $\k_f$-linear isomorphisms, 
$\lambda_f\colon f^k \L/f^{k+1}\L \isom f^{k+1} \L/f^{k+2}\L$,
for all $k\ge 0$.
\end{proof}

\begin{corollary}
\label{cor:iso}
The map $\lambda_f\colon \wL \to \wL$ is injective 
(that is, $f$ is not a zero-divisor in $\wL$).
\end{corollary}

\begin{proof}
Due to the above Lemma, $\wL$ is an integral domain; see
\cite[Proposition 7 on p.~II-8]{Se75}. Our claim follows then 
from the injectivity of $\iota$.
\end{proof}

\subsection{$\Z$-covers}
\label{subsec:zcov}

Now let $X$ be a connected CW-complex, and let 
$\nu\colon \pi_1(X)\to \Z$ be a non-trivial homomorphism.    
Since $\Z$ is abelian, $\nu$ factors through $H_1(X,\Z)$, 
thus defining a $1$-dimensional cohomology class, 
$\nu_{\Z}\in H^1(X,\Z)$.  Let $\omega\in H^1(S^1,\Z)=\Z$ 
be the standard generator.  By obstruction theory, there 
is a map $h\colon X\to S^1$ such that $\nu_{\Z}=h^*(\omega)$.  
It is readily seen that the induced homomorphism, 
$h_{\sharp} \colon \pi_1(X)\to \pi_1(S^1)$, 
coincides with $\nu$.   

Denote by $X^{\nu}$ the Galois (connected) cover of $X$ 
corresponding to $\ker (\nu)$.  Without loss of generality, 
we may assume $\nu$ is surjective.  Indeed, if the image 
of $\nu$ has index $m$ in $\Z$, then $\nu$ is the composite 
$\pi_1(X)\xrightarrow{\nu'} \Z \xrightarrow{\times m} \Z$, 
where $\nu'$ is onto.  Clearly, $\ker(\nu)=\ker(\nu')$, and 
so the corresponding covers are equivalent. 

We may obtain the  infinite cyclic Galois cover,  
$\pi\colon X^{\nu} \to  X$, with $\im(\pi_{\sharp})=\ker(\nu)$, by
pulling back the universal cover $\exp\colon \R\to S^1$ 
along $h$. 
Note that $X^{\nu}$ is the homotopy fiber of $h$. The 
spaces and maps defined so far fit into the following 
diagram:
\begin{equation}
\label{eq:cd}
\xymatrix{
& \Z \ar[d] \ar[r]& \Z \ar[d] \\
X^{\nu} \ar@{=}[d] \ar@{=}[r] 
& X^{\nu} \ar^{\pi}[d]\ar[r] & \R \ar^{\exp}[d] \\
X^{\nu} \ar^{\pi}[r] & X \ar^{h}[r] & S^1 \\
}
\end{equation}

\subsection{Cohomology of $\Z$-covers}
\label{subsec:coho zcover}

Let $\k$ be a coefficient field.  The cell 
structure on $X$ lifts to an equivariant cell structure on 
$X^{\nu}$, so that the group $\Z$ acts by deck transformations 
on $X^{\nu}$, permuting the cells.  In this fashion, each 
homology group $H_{i}(X^{\nu},\k)$ acquires the structure 
of a module over the group ring $\k\Z$.  Since $\nu$ 
is surjective, Shapiro's lemma yields an isomorphism of 
$\k\Z$-modules, 
\begin{equation}
\label{eq:shapiro}
H_{i}(X^{\nu} ,\k) \cong  H_{i}(X ,\k\Z_{\nu}),
\end{equation}
where $\k\Z_{\nu}$ denotes the ring $\k\Z$, viewed 
as a module over $\k\pi_1(X)$ via the linear extension 
$\tilde\nu\colon \k\pi_1(X)\to \k\Z$. 

Let $\nu_{*}\colon H_1(X,\k)\to H_1(\Z,\k)=\k$ be the 
homomorphism induced by $\nu$. The corresponding 
cohomology class,  $\nu_{\k}\in H^1(X,\k)$, is the image 
of $\nu_{\Z}$ under the coefficient homomorphism $\Z\to \k$. 
Since $\nu_{\Z}\cup\nu_{\Z}= h^*(\omega\cup\omega) =0$, 
we also have $\nu_{\k} \cup \nu_{\k} = 0$, by naturality of cup 
products with respect to coefficient homomorphisms. 
We will denote by $(\nu_\k)$ the ideal of $A=H^*(X,\k)$ 
generated by $\nu_\k$.  Moreover, we will write 
$A^{\le r}$ for the ring $A$, modulo the ideal 
$\bigoplus_{i>r} A^i$; of course, additively 
$A^{\le r}=\bigoplus_{i=0}^r A^i$. 

\begin{prop}
\label{prop:coho zcover}
Let $\pi\colon X^{\nu} \to X$ be a $\Z$-cover as above. 
Then:
\begin{enumerate}
\item \label{coho1}
The induced homomorphism in cohomology, 
$\pi^{*}\colon H^{*}(X,\k) \to H^{*}(X^{\nu},\k)$, 
factors through a ring map, 
$\bar\pi^*\colon H^{*}(X,\k)/(\nu_\k) \longrightarrow  
H^{*}(X^{\nu},\k)$.
\item \label{coho2}
Suppose $H_{i}(X^{\nu},\k)$ has trivial $\Z$-action, 
for all $i\le r$.  Then $\bar\pi^*$ restricts to a ring isomorphism, 
$\bar\pi^*\colon H^{\le r}(X,\k)/(\nu_\k) \isom 
H^{\le r}(X^{\nu},\k)$.
\end{enumerate}
\end{prop}

\begin{proof}
Consider the cohomology spectral sequence with 
$\k$-coefficients associated to the homotopy fibration 
$X^{\nu}\xrightarrow{\pi} X \xrightarrow{h}  S^1$. 
Clearly, $E_2=E_{\infty}$. Note also that the $\Z$-action 
on the homology of $X^{\nu}$ given by \eqref{eq:shapiro} 
may be identified with the action of $\pi_1(S^1)$ on the 
homology of the fiber, and similarly for cohomology. 
The factorization in Part \eqref{coho1} follows from the fact that
$\nu_{\k}= h^*(\omega_{\k})$.
The claim in Part \eqref{coho2} now follows from \cite{Bo}: 
the surjectivity of $\bar{\pi}^*$ from (13.5) on p.~39, 
and the injectivity of $\bar{\pi}^*$ from Theorem 14.2(b) 
on p.~42.
\end{proof}

\subsection{Field of fractions}
\label{subsec:fractions}

For the rest of this section, we shall assume $X$ has 
finite $k$-skeleton, for some $k\ge 1$; in particular, 
$\pi_1(X)$ is finitely generated.  We seek to compute 
the finitely generated $\k\Z$-modules $H_{i}(X^{\nu},\k)$, 
for $i\le k$. As a first step, we reduce the computation 
of the ranks of these $\k\Z$-modules to that of certain 
Betti numbers with coefficients in rank $1$ local systems. 

Let $\KK=\k(t)$ be the field of fractions of the integral domain 
$\k\Z=\k[t^{\pm 1}]$. Set $G:= \pi_1(X)$. Composing the 
natural inclusion $\k\Z \inj \KK$ with the linear extension 
of $\nu$ to group rings, $\k{G}\to \k\Z$, and restricting to 
$G$ defines a homomorphism $\bar\nu\colon G \to \KK^{\times}$.  
In this manner, we have a rank $1$ local system, 
$\KK_{\bar\nu}$, on $X$. 

\begin{lemma}
\label{lem:hom nu}
With notation as above,
\[
\rank_{\k\Z} H_{i}(X^{\nu} ,\k) = 
\dim_{\KK} H_{i}(X, \KK_{\bar\nu}),
\quad \text{for $i\le k$} .
\]
\end{lemma}

\begin{proof}
Note that the inclusion $\k\Z \inj \KK$ is a flat morphism; 
thus,  $\rank_{\k\Z} H_{i}(X ,\k\Z_{\nu}) = 
\dim_{\KK} H_{i}(X,\k\Z_{\nu} \otimes_{\k\Z} \KK)$. 
Since $\k\Z_{\nu} \otimes_{\k\Z} \KK=\KK_{\bar\nu}$, 
we are done. 
\end{proof}

\subsection{Completion}
\label{subsec:comp}

In each degree $0\le i\le k$, the finitely generated 
$\L$-module $H_i(X^{\nu},\k)$ decomposes as
\begin{equation}
\label{eq:hq decomp}
H_i(X^{\nu} ,\k) =\L^{r_i} \oplus 
\bigoplus_{\stackrel{t\ne f\in \k[t]}{f\:\text{irreducible}}} 
\bigoplus_{j\ge 1} \left(\L/f^j \L\right)^{e^i_j(f)},
\end{equation}
where $r_i$ is the $\L$-rank, and $e^i_j(f)$ is 
the number of $f$-primary Jordan blocks of size $j$. 
Of course, since $X^{\nu}$ is path-connected, 
$H_0(X^{\nu},\k)=\L/(t-1)$, and so $r_0=0$, and 
$e^0_1 (t-1)=1$ is the only non-zero entry among 
the multiplicities $e^0_j(f)$. 

Now fix an irreducible polynomial $f\in \L$, and let 
$\wL$ be the completion of $\L$ with respect to the 
$(f)$-adic filtration. Recall we have an injective map 
$\iota\colon \L\to \wL$. Since $\iota$ is a flat morphism, 
decomposition \eqref{eq:hq decomp} yields, for each $i\le k$, 
\begin{align}
\label{eq:hom complete}
H_i(X,\wL_{\nu}) 
&= H_i(X,\L_{\nu}) \otimes_{\L} \wL \\
&= \wL^{r_i} \oplus \bigoplus_{j\ge 1} 
\big(\wL/f^j \wL\big)^{e^i_j(f)}.\notag
\end{align}

Note that this gives the canonical decomposition of 
$H_i(X,\wL_{\nu})$ over the principal ideal domain 
$\wL$, since $f$ is (the unique) irreducible in $\wL$.

\section{Homology of $\Z$-covers of toric complexes}
\label{sect:mono}

We now specialize to the case when $X$ is a toric 
complex, and compute the homology groups of an infinite 
cyclic cover, viewed as modules over the ring of Laurent 
polynomials.  

\subsection{$\Z$-covers and Artin kernels}
\label{subsec:zcovers}

As before, let $L$ be a finite simplicial complex, 
on vertex set $\sV$, and let $T_L$ be the corresponding 
toric complex. Recall that the fundamental group
$G_L=\pi_1(T_L)$ depends only on the graph $\G=L^{(1)}$, 
and that a classifying space for $G_L$ is the toric complex 
$T_{\Delta}$, where $\Delta=\Delta_{\G}$ is the flag complex 
of $\G$. 

Consider now an epimorphism $\chi\colon G_{L}\surj \Z$. 
The construction reviewed in \S\ref{subsec:zcov} defines 
an infinite cyclic (regular) cover, $\pi\colon T_L^{\chi} \to T_{L}$. 
The fundamental group 
\begin{equation}
\label{eq:ak def}
N_{\chi}:=\pi_1(T_L^{\chi})=\ker (\chi\colon G_{L}\to \Z)
\end{equation}
is called the {\em Artin kernel}\/ associated to $\chi$.  
Clearly, a classifying space for this group is the space 
$T_{\Delta}^{\chi}$. 

The most basic example is provided by the ``diagonal" 
homomorphism, $\nu\colon G_L\surj \Z$, given by 
$\nu(v)=1$, for all $v\in \V$.  The corresponding 
Artin kernel, $N_{\nu}$, is simply denoted by $N_\G$, 
and is called the {\em Bestvina-Brady group} associated 
to $\G$.   This group need not be finitely generated. 
For example, if $\overline{K}_n$ is the discrete graph 
on $n>1$ vertices, then $G_{\overline{K}_n}=F_n$ 
and $N_{\overline{K}_n}$ is a free group of 
countably infinite rank.  More generally, as shown by 
Meier--VanWyk \cite{MV} and  Bestvina--Brady \cite{BB} 
the group $N_{\G}$ is finitely generated if and only 
if the graph $\G$ is connected. Even then, $N_{\G}$ 
may not admit a finite presentation. For example, if 
$K_{2,2}=\overline{K}_2* \overline{K}_2$ is a $4$-cycle, 
then $G_{K_{2,2}}= F_2\times F_2$, 
and, as noted by Stallings \cite{St}, 
$N_{K_{2,2}}$ is not finitely presentable. In fact, as 
shown in \cite{BB}, $N_{\G}$ is finitely presented if 
and only if $\Delta_\G$ is simply-connected

More generally, we have the following characterization 
from Meier--Meinert--VanWyk \cite{MMV} and 
Bux--Gonzalez \cite{BG}. Assume $L$ is a flag complex. 
Let $\sW= \{ v\in \sV \mid \chi(v) \ne 0 \}$ 
be the support of $\chi$. Then:
\begin{alphenum}
\item \label{test1} 
$N_{\chi}$ is finitely generated if and only if 
$L_{\sW}$ is connected, and dominant, i.e., for all 
$v\in \sV\setminus \sW$, there is a $w\in \sW$ such 
that $\set{v,w}\in L$.   
\item \label{test2} 
$N_{\chi}$ is finitely presented if and only if 
$L_{\sW}$ is $1$-connected and, for every 
simplex $\sigma$ in $L_{\sV\setminus \sW}$, the space  
$\lk_{L_{\sW}}(\sigma):= \{ \tau\in L_{\sW} 
\mid \tau \cup \sigma \in L \}$ 
is $(1-\abs{\sigma})$-acyclic. 
\end{alphenum}

\subsection{$\k \Z$-ranks for $\Z$-covers of toric complexes}
\label{subsec:zcov tc}
Let $L$ be an arbitrary finite simplicial complex.
Fix a cover $T_L^{\chi}\to T_L$, defined by a homomorphism 
$\chi\colon G_L\surj \Z$. Our goal in this section is to compute 
the homology groups $H_i(T_{L}^{\chi} ,\k)$, viewed as modules 
over the group ring $\k\Z$, for a fixed coefficient field $\k$. 
In view of  formula \eqref{eq:hq decomp}, we need to compute 
the $\k\Z$-ranks, $r_i$, and the multiplicities of the Jordan 
blocks, $e^i_j(f)$, for all integers $i\ge 0$, $j\ge 1$, and 
all irreducible polynomials $f\in \k[t]$, $f\ne t$. We begin 
with a combinatorial formula for the $\k\Z$-ranks.

Define $\sV_{\k} (\chi) =\supp( \chik)$, where 
$\chik\in H^1(T_L,\k)=\SR^{1}$ is the corresponding 
cohomology class.  Clearly, the subset $\sV_{\k}(\chi)\subset \sV$ 
depends only on $\chi$ and $p=\ch\k$, so we simply write 
\begin{equation}
\label{eq:supp nu}
\sV_{p} (\chi):=\sV_{\k} (\chi) =\supp( \chik). 
\end{equation}

\begin{theorem}
\label{thm:hom cov}
For all $i\ge 0$,
\[
\rank_{\k\Z} H_i(T^{\chi}_L ,\k) = \beta_i (\SR ,\sV_0(\chi)).
\]
\end{theorem}

\begin{proof}
Set $\KK=\k(t)$, and let $\bar\chi\colon G_{L} \to \KK^{\times}$ 
be the $1$-dimensional representation defined by $\chi$. 
The cohomology class 
$z_{\bar\chi}=\sum_{v\in \sV}(\bar\chi(v)-1) v^* \in  
\KK\langle L \rangle^{1}$ has support
\begin{equation}
\label{eq:supps}
\supp (z_{\bar\chi}) =\set{ v\in \sV \mid \bar\chi(v)\ne 1}
=\set{ v\in \sV \mid \chi(v)\ne 0} = \sV_0(\chi). 
\end{equation}

We have the chain of equalities 
\begin{align*}
\label{eq:hom eq}
\rank_{\k\Z} H_i(T^{\chi}_L ,\k) 
&= \dim_{\KK} H_i(T_L, \KK_{\bar\chi}) 
&\text{by Lemma \ref{lem:hom nu}}\\
&=\beta_i (\KK\langle L \rangle , z_{\bar\chi}) 
&\text{by Lemma \ref{lem:change vars}} \\
&= \beta_i (\KK\langle L \rangle ,\supp(z_{\bar\chi}) )
&\text{by Lemma \ref{lem:aah}}\\
&=\beta_i (\KK\langle L \rangle , \sV_0(\chi)) 
&\text{by \eqref{eq:supps}}\\
&= \beta_i (\SR ,\sV_0(\chi)) 
&\text{since $\ch\KK=\ch\k$}.
\end{align*}
This ends the proof.
\end{proof}

Here is an immediate consequence of Theorem \ref{thm:hom cov} 
(together with Corollary \ref{cor:beta vanish} and 
Remark \ref{rem:beta vanish}). 
\begin{corollary}
\label{cor:hom res}
For each $r\ge 1$, the following are equivalent:
\begin{enumerate}
\item \label{bns1}  $\dim_{\k} H_{i} (T_L^{\chi} ,\k) < \infty$, 
for all $i\le r$. 
\item \label{bns2}  $\beta_i (\SR,\sV_0 (\chi))=0$, for all $i\le r$.
\item \label{bns3}  
$\widetilde{H}_{i} (\lk_{L_{\sV_0(\chi)}}(\sigma),\k)=0$, 
for  all $\sigma\in L_{\sV\setminus \sV_0(\chi)}$ and 
$-1\le i\le r-1-\abs{\sigma}$.
\end{enumerate}
If $\ch\k=0$, the above conditions are also equivalent to
\begin{enumerate}
\setcounter{enumi}{3}
\item \label{bns4}  $\chik\notin \RR^i_1(T_L,\k)$, for all $i\le r$. 
\end{enumerate}
\end{corollary}

\subsection{A chain isomorphism}
\label{subsec:chain iso}

Let $R$ be an arbitrary commutative ring, and suppose we 
are given a vector $\gamma =(\gamma_v )_{v\in \V}\in R^{\V}$.
Let $C_{\bullet}(T_L, R_{\gamma})$ be the $R$-chain complex 
defined as follows. Set
$C_k(T_L,R_{\gamma})=R\otimes C_k$, where 
$C_k=C_k(T_L,\Z)$ are the usual cellular $k$-chains on $T_L$.  
The boundary maps $\partial^{\gamma}\colon R\otimes C_k \to 
R\otimes C_{k-1}$ are given by 
\begin{equation}
\label{eq:bdry map}
\partial^{\gamma} (1\otimes c_{\sigma})=  
\sum_{r=1}^{k} (-1)^{r-1} \gamma_{v_r}
\otimes c_{\sigma\setminus \{v_r\}}, 
\end{equation}
where $\sigma=\{v_{1},\dots , v_{k}\}$ is a $(k-1)$-simplex 
in $L$. 

Given two vectors $\gamma=(\gamma_v)_{v\in \V}$ 
and  $\xi=(\xi_v)_{v\in \V}$ in $R^{\V}$, let 
us write $\gamma \doteq \xi$ if there exist units 
$u_v\in R^{\times}$ such that $\gamma_v=u_v \xi_v$, 
for all $v\in \V$. 

\begin{lemma}
\label{lem:chain iso}
If $\gamma \doteq \xi$, then $H_*(T_L, R_{\gamma}) 
\cong H_*(T_L, R_{\xi})$, as $R$-modules.  
\end{lemma}

\begin{proof}
Define a chain map $C_{\bullet}(T_L, R_{\gamma})\to 
C_{\bullet}(T_L, R_{\xi})$ by $1\otimes c_{\sigma} 
\mapsto u_{v_1}\cdots u_{v_k} \otimes c_{\sigma}$, for 
$c_{\sigma}\in C_k$.  Clearly, this is a chain isomorphism.  
\end{proof}

The computation of the $\k \Z$-torsion part of $H_{i} (T_L^{\chi} ,\k)$,
encoded by the multiplicities $e^i_j(f)$, is in two steps: 
the first, arithmetic, and the second, combinatorial.

\subsection{The arithmetic step}
\label{subsec:d at f}

Fix an irreducible polynomial $f\in \k[t]$, different from $t$. 
Set $m_v =\chi (v)\in \Z$, for each $v\in \V$.  Define the 
vector $b=(b_v)_{v\in \V}$, with components 
\begin{equation}
\label{eq:bv}
b_v=b_v(\chi,f)= \begin{cases}
\ord_f (t^{m_v} -1), & \text{if $m_v\ne 0$},\\
-\infty, & \text{if $m_v= 0$},
\end{cases}
\end{equation}
where $\ord_f(g)=\max\{ k \ge 0: f^k \!\mid g\}$, for 
$0\ne g\in \L$.  

\begin{example}
\label{ex:alg closed}
To illustrate the computation, suppose $\k$ is an algebraically 
closed field, of characteristic $p$. Then any irreducible 
polynomial in $\L$ is of the form $t-a$, for some 
$a\in \k^{\times}$. Assume $m_v\ne 0$. If either 
$\ord(a)=\infty$, or $\ord(a)<\infty$ and $\ord(a)\nmid m_v$, 
then $b_v=0$; otherwise, 
\begin{equation}
\label{eq:bv formula}
b_{v} = \begin{cases}
1, & \text{if $p=0$},\\
p^s, & \text{if $m_v=p^sq$, with $(q,p)=1$}.
\end{cases}
\end{equation}
\end{example}

Returning now to the general case, denote by 
$\wL$ the $(f)$-adic completion of $\L$.  Define 
a new vector, $\gamma =\gamma (\chi, f)\in \wL^{\V}$, 
by setting $\gamma_v =f^{b_v}$, with the convention 
that $f^{-\infty}=0$. The associated differential, 
$\partial^{\gamma}\colon \wL \otimes C_{k}\to 
\wL \otimes C_{k-1}$, is given by 
\begin{equation}
\label{eq:dbf}
\partial^{\gamma} (1\otimes c_{\sigma})=  
\sum_{r=1}^{k} (-1)^{r-1} f^{b_{v_r}} 
\otimes c_{\sigma\setminus \{v_r\}}. 
\end{equation}
Using formula \eqref{eq:equiv bdry} and Remark 
\ref{rem:units}, we see that $\partial^{\gamma}$ equals 
$\id_{\wL_{\chi}}\otimes_{\k{G_L}} \widetilde{\partial}$, 
modulo units in $\wL$.  From Lemma \ref{lem:chain iso}, 
we derive the following corollary. 

\begin{corollary}
\label{cor:hom bv}
$H_*(T_L, \wL_{\chi})= H_*(T_L, \wL_{\gamma})$, as $\wL$-modules.
\end{corollary}

\subsection{Independence of $f$}
\label{subsec:indep}

The output of the preceding step is the vector 
$b=(b_v)_{v\in \sV}$, constructed in \eqref{eq:bv} via 
the factorization properties of the ring $\L$.  We will use this 
vector as input in the second, combinatorial, step of our 
$\k\Z$-torsion computation, in the following way.

Let $S=\k[[t]]$ be the power-series ring in variable $t$.
Define a vector $\xi=(\xi_v)_{v\in \V}\in S^{\V}$ by 
$\xi_v =t^{b_v}$, where $t^{- \infty}:= 0$. Since 
$S$ is a PID with a unique prime element, $t$, the 
$S$-module $H_{i}(T_L, S_{\xi})$ decomposes as 
\begin{equation}
\label{eq:s decomp}
H_{i}(T_L, S_{\xi})= S^{\rho_{i}} \oplus \bigoplus_{j\ge 1} 
\big(S/t^j S\big)^{\varepsilon^i_j} ,
\end{equation}
for all $i\ge 0$. The next result shows, in particular, 
that the multiplicities $e^i_j(f)$ do not depend on $f$.

\begin{prop}
\label{prop:ttors}
With notation as above, $r_i=\rho_i$ and 
$e^i_j(f)=\varepsilon^i_j$, for all $i\ge 0$ and $j\ge 1$.
\end{prop}

\begin{proof}
Since $S$ is the $(t)$-adic completion of $\k [t]$, we have 
a canonical ring map $S\to \wL$, $t\mapsto f$. It follows 
from \eqref{eq:bdry map} that
\begin{equation}
\label{eq:extend}
C_{\bullet}(T_L, \wL_{\gamma}) = 
\wL \otimes_S C_{\bullet}(T_L, S_{\xi}).
\end{equation}
Using the Universal Coefficients Theorem (over $S$) 
and decomposition \eqref{eq:s decomp}, we find 
$\wL$-isomorphisms
\begin{align}
\label{eq:hqs}
H_{i}(T_L,\wL_{\gamma}) 
&= \big( \wL \otimes_S H_{i}(T_L, S_{\xi})\big) \oplus 
\Tor_1^S (\wL, H_{i-1}(T_L, S_{\xi}))\\
&= \wL^{\rho_{i}} \oplus  
\Big( \bigoplus_{j\ge 1} (\wL/f^j \wL)^{\varepsilon^i_j}\Big)
\oplus \Big( \bigoplus_{j\ge 1} \ker (\lambda_{f^j}\colon 
\wL \to \wL)^{\varepsilon^{i-1}_j}\Big)\, .
\notag
\end{align}

By Corollary \ref{cor:iso}, the third summand above vanishes. 
Using Corollary \ref{cor:hom bv}, our claim follows by  
comparing the $\wL$-decompositions \eqref{eq:hqs} 
and \eqref{eq:hom complete}.
\end{proof}

\subsection{Realizability}
\label{subsec:discuss}
Before proceeding to the combinatorial step of our 
algorithm, we discuss the possible torsion that can 
occur in the $\k\Z$-module $H_*(T_L, \k\Z_{\chi})$. 

\begin{prop}
\label{prop:triv fprimary}
Let $L$ be a simplicial complex, and let $f$ be a 
polynomial in $\k[t]$, irreducible and different from $t$.  
Suppose that, for every vertex $v\in \sV$ such that $m_v\ne 0$, 
the polynomial $f$ does not divide $t^{m_v}-1$. Then 
$H_*(T_L, \k\Z_{\chi})$ has trivial $f$-primary part.  
\end{prop}

\begin{proof}
Let $(b_v)_{v\in \sV}$ be the vector defined by \eqref{eq:bv}. 
Our assumption implies that $b_v=0$, if $m_v\ne 0$. 
Define a $\k$-valued vector $(\delta_v)_{v\in \V}$ by 
$\delta_v=1$, if $m_v\ne 0$, and $\delta_v=0$,
otherwise. Let $\gamma$ and $\xi$ be the vectors 
defined in \S\ref{subsec:d at f} and \S\ref{subsec:indep}. 
Using formula \eqref{eq:bdry map} for the boundary 
map $\partial^{\xi}$, we see that 
$H_i(T_L, S_{\xi})$ is a free $S$-module of rank 
equal to $\dim_{\k} H_i(T_L, \k_{\delta})$, for all $i\ge 0$. 
The claim follows from Proposition \ref{prop:ttors}. 
\end{proof}

Conversely, we have the following. 

\begin{prop}
\label{prop:realize}
Let $f$ be an irreducible polynomial in $\k[t]$ dividing 
$t^m -1$, for some $m\ge 1$.   Then, for any $i>0$, there 
exists a simplicial complex $L$ and a homomorphism 
$\chi\colon G_{L}\surj \Z$ such that $H_i(T_L, \k\Z_{\chi})$ 
has nontrivial $f$-primary part. 
\end{prop}

\begin{proof}
Let $L$ be the cone $v_0*K$, where $K =\Delta_{\Gamma}$ 
is a flag triangulation of $S^{i-1}$.  Define $\chi (v)=1$ on the 
vertices of $K$, and $\chi(v_0)=m$. By the K\"{u}nneth formula, 
the $\k\Z$-module $H_i(T_L, \k\Z_{\chi})$ has 
$M=\k\Z/(t^m-1)\otimes_{\k\Z} H_i(N_{\Gamma}, \k)$ as a 
direct summand. By  \cite[Proposition 7.1]{PS-bb}, the module 
$M$ has $\k\Z/(t^m-1)$ as a direct $\k\Z$-summand. 
Clearly, $\k\Z/(t^m-1)$ has non-trivial $f$-primary part, 
and we are done.
\end{proof}

\subsection{The combinatorial step}
\label{subsec:comb step}
The computation of $T_f(H_{i} (T_L,\k\Z_{\chi}))$ 
described in \S\ref{subsec:indep} depends only on the 
simplicial complex $L$, and on the vector $b=(b_v(\chi,f))$.  
More precisely, given the vector $\xi=(t^{b_v})_{v\in\sV} 
\in S^{\V}$ as above, the torsion part of the $S$-decomposition 
\eqref{eq:s decomp} can be computed from the boundary 
map $\partial_{i+1}^{\xi}\colon 
S\otimes C_{i+1}\to S\otimes C_{i}$, given by 
\begin{equation}
\label{eq:dt}
\partial^{\xi} (1\otimes c_{\sigma})=  
\sum_{r=1}^{i+1} (-1)^{r-1} t^{b_{v_r}} 
\otimes c_{\sigma\setminus \{v_r\}}. 
\end{equation}

Indeed, consider the split exact sequence of $S$-modules, 
\begin{equation}
\label{eq:les}
\xymatrix{0\ar[r]& \ker(\partial_{i}^{\xi})/\im (\partial_{i+1}^{\xi}) \ar[r]& 
\coker(\partial_{i+1}^{\xi})\ar[r]& \im 
(\partial_{i}^{\xi}) \ar[r]& 0}.
\end{equation}
It is readily seen that $\Tors H_i(T_L, S_{\xi})= 
\Tors \coker(\partial_{i+1}^{\xi})$.  Since $S$ 
is a PID, the matrix of $\partial_{i+1}^{\xi}$ 
can be diagonalized, using elementary row and 
column operations. In other words, with respect to 
convenient bases over $S$, the matrix of 
$\partial_{i+1}^{\xi}$ can be written as 
$0\oplus D_{\xi}$, where
\begin{equation}
\label{eq:diag matrix}
D_{\xi}=\diag (t^{a_k}), \quad \text{with $a_k \ge 0$}.
\end{equation}
Clearly, 
\begin{equation}
\label{eq:multiplicities}
\varepsilon^i_j=\abs{\set{k \mid a_k= j}}.
\end{equation}
Using Proposition \ref{prop:ttors}, we conclude that 
$T_f(H_{i} (T_L,\k\Z_{\chi}))= \bigoplus_{k} (\L/f^{a_k} \L)$. 

The next Theorem summarizes our two-step algorithm 
for computing the $\k\Z$-torsion in $H_{i}(T_L^{\chi},\k)$, 
for all $i\ge 0$.

\begin{theorem}
\label{thm:hom tors}
Fix an irreducible polynomial $f\in \k[t]$, $f\ne t$. To compute 
the $f$-primary part $T_f(H_{i} (T_L^{\chi},\k))=\bigoplus_{j\ge 1} 
(\Lambda/f^j \Lambda)^{e^i_j(f)}$, proceed as follows. 
\begin{enumerate}
\item \label{step1}
Compute the vector $b=(b_v)_{v\in \V}$, with $b_v=b_v(\chi,f)$ 
given by \eqref{eq:bv}. 
\item \label{step2}
With this vector as input, define the boundary map 
$\partial^{\xi}$ as in \eqref{eq:dt}, and compute 
the multiplicities $\varepsilon^i_j$ as in \eqref{eq:multiplicities}, 
using the diagonalized matrix $D_{\xi}$ from 
\eqref{eq:diag matrix}.
\end{enumerate}
Finally, set $e^i_j(f)=\varepsilon^i_j$. 
\end{theorem}

\subsection{The Bestvina-Brady covers}
\label{subsec:bb covers}
Particularly simple is the case when $\chi$ is the 
diagonal homomorphism $\nu\colon G_L\surj \Z$, 
taking each $v\in \sV$ to $1$. The $\k\Z$-module structure 
on the homology of the resulting Bestvina-Brady cover, 
$T_L^{\nu}$, can be made completely explicit.  

\begin{corollary}
\label{cor=bbcover}
For each $i>0$, there is an isomorphism of $\L$-modules,
\[
H_{i} (T_L^{\nu},\k)= 
\L^{\dim_{\k} \widetilde{H}_{i-1}(L, \k)} \oplus
\big( \L/(t-1)\L \big)^{\dim_{\k} B_{i-1}(L, \k)} ,
\]
where $B_{\bullet}(L, \k)$ are the simplicial 
boundaries of $L$.
\end{corollary}

\begin{proof}
By Theorem \ref{thm:hom cov} the free part of 
$H_{i} (T_L^{\nu},\k)$ has rank $\beta_i(\SR, \V)$.   
By Proposition \ref{prop:hochster aah}, this equals 
$\dim_{\k} \widetilde{H}_{i-1}(L, \k)$. 

By Proposition \ref{prop:triv fprimary}, the torsion part 
is all $(t-1)$-primary.  The arithmetic step \eqref{eq:bv} 
is very simple: it gives $\xi_v=t$, for all $v\in \V$.  
Therefore, the diagonal matrix $D_{\xi}$ has size 
$\dim_{\k} B_{i-1}(L, \k)$, with all diagonal entries 
equal to $t$.  The conclusion follows.
\end{proof}

When applied to a flag complex $L=\Delta_\G$, 
Corollary \ref{cor=bbcover} recovers Proposition 7.1 
from \cite{PS-bb}.

\section{Trivial monodromy test}
\label{sect:tm test}

In this section, we give a combinatorial test for 
deciding whether $H_{\le r} (T_L,\k\Z_{\chi}):= \bigoplus_{i=0}^{r} H_{i} (T_L,\k\Z_{\chi})$ 
has trivial $\k\Z$-action.

\subsection{Supports and primes}
\label{subsec:supp}
First, we need to establish some notation.  Let 
$\chi\colon G_L\surj \Z$ be an epimorphism. 
This amounts to specifying integers $\chi(v)=m_v$ 
for each vertex $v\in \sV$, with the proviso that 
$\gcd\set{m_v\mid v\in \sV}=1$. In other words, 
the epimorphism $\chi\colon G_L\surj \Z$ 
is encoded by a graph $\G=L^{(1)}$, 
equipped with a vertex-labeling function 
$m\colon \sV \to \Z$, subject to the 
above coprimality condition. A simple example 
of such a vertex-labeled graph is given in diagram 
\eqref{eq:path}. 

In this setup, the support sets from \eqref{eq:supp nu} 
can be written as
\begin{equation}
\label{eq:vq}
\sV_q(\chi)=\set{v \in \sV\mid m_v \ne 0 \!\pmod{q}}, 
\end{equation}
for $q=0$ or $q$ a prime. 
Clearly, $\sV_q(\chi)\subset \sV_0(\chi)\subset \sV$. 
Furthermore, $\sV_q(\chi)\ne \emptyset$, since the 
${m_v}'s$ are coprime.  Also define
\begin{equation}
\label{eq:pnu}
\PP(\chi)=\set{\text{$q$ prime}\mid \sV_q(\chi)
\subsetneq \sV_0(\chi)}.
\end{equation}
Obviously, this is a finite (possibly empty) set, consisting 
of all the prime factors of the non-zero ${m_v}'s$. 

Now let $\k$ be a field, of characteristic $p$. 
For an element $a\in \k\setminus \set{0,1}$, 
the homomorphism 
$\chi_a:=\ev_a \circ \chi \colon G_L\to \k^{\times}$ 
is given by $\chi_a(v)=a^{m_v}$.  Let 
\begin{equation}
\label{eq:znua}
z_{\chi_a}=\sum_{v\in \sV}(a^{m_v}-1) v^* \in  \SR^{1}
\end{equation}
be the corresponding cohomology class. Set 
$\sW_a=\supp(z_{\chi_a})$. Clearly,  $v\in \sW_a$ 
if and only if $a^{m_v}\ne 1$ in $\k$. Thus,
\begin{equation}
\label{eq:wa}
\sW_a=\{v\in \sV_0(\chi)\mid \ord(a) \nmid m_v\}.
\end{equation}

\begin{lemma}
\label{lem:supp}
With notation as above, 
\begin{enumerate}
\item \label{s1}
If $\ord(a)=\infty$, then $\sW_a=\sV_0(\chi)$.
\item \label{s2}
If $\ord(a)<\infty$, then there is a prime $q\ne p$ 
such that $\sW_a\supset \sV_q(\chi)$.
\end{enumerate}
\end{lemma}

\begin{proof}
If $a$ has infinite order, then clearly 
$\sW_a= \set{v\in \sV\mid m_v\ne 0}$, which equals 
$\sV_0(\chi)$. If $a$ has finite order $d$ (necessarily, $d>1$, since 
$a\ne 1$), then $\sW_a=\set{v\in \sV_0 (\chi)\mid d \nmid m_v}$. 
We claim there is a prime $q\ne p$ such that $q\mid d$.  
If $p=0$, this is clear.  If $p>0$ and $d=p^s$, 
then $0=a^d-1=(a-1)^{p^s}$, which forces $a=1$, 
a contradiction. Hence, 
$\sV_q(\chi)=\set{v\in \sV_0(\chi)\mid q\nmid m_v}
\subset \sW_a$. 
\end{proof}

\subsection{Triviality test for the monodromy}
\label{subsec:triv mono}
We are now ready to state and prove the main 
result of this section.

\begin{theorem}
\label{thm:mono test}
Let $\k$ be a field, with $\ch \k=p$. 
The $\k\Z$-module $H_{\le r} (T_L,\k\Z_{\chi})$ 
has trivial $\Z$-action if and only if 
\begin{enumerate}
\item[$(\dagger)_r$] 
$\beta_i(\SR,\sV_p(\chi))=0$, for all $i\le r$;
\item[$(\ddagger)_r$] 
$\beta_i(\SR,\sV_q(\chi))=0$, for all $q\in \PP(\chi)$ 
with $q\ne p$, and for all $i\le r$. 
\end{enumerate}
\end{theorem}

\begin{proof}
Everything depends only on $p=\ch \k$, and not on $\k$ itself, 
so we may as well assume $\k$ is algebraically closed. 
It follows that, up to units, the only irreducible polynomials in 
$\Lambda=\k\Z$ are of the form $t-a$, with $a\in \k^{\times}$.  
Thus, we may write 
\begin{equation}
\label{eq:afact}
H_i(T_L,\L_{\chi}) =\L^{r_i} \oplus T^i_1 \oplus 
\bigoplus_{a\ne 1} T^i_a,
\end{equation}
where $T_{a}$ denotes $(t-a)$-primary part. 

Let  $\chik\in\SR^1=H^1(T_L,\k)$ be the cohomology class 
corresponding to the composite $G_L\xrightarrow{\chi} \Z \to \k$. 
Clearly, $\supp(\chik)= \V_p(\chi)$; thus, 
$\beta_i(\SR,\chik)= \beta_i(\SR, \V_p(\chi))$,
by Lemma \ref{lem:aah}.  It follows from 
\cite[Proposition 9.4]{PS-eqchain} that 
\begin{align}
\label{eq:mono eqchain}
\text{$\L^{r_i} \oplus T^i_1$ has trivial
$\Z$-action, $\forall i\le r$} &\same 
\text{$r_i=0$ and $T^i_1=(\Lambda/(t-1))^{e_i}$, 
$\forall i\le r$}\\
&\same \text{condition $(\dagger)_r$ is satisfied}.\notag
\end{align}

Now fix $a\ne 0$ or $1$.  Since $t-1$ and $t-a$ are coprime, 
$t-1$ acts trivially on $T_a$ if and only if $T_a=0$.  
So we're left with proving the following two claims.

\begin{claim}
\label{claim:1}
Suppose condition $(\dagger)_{r}$ is satisfied. 
If $T^{\le r}_a=0$, for all $a\in\k\setminus\set{0,1}$,  
then $\beta_{i}(\SR,\sV_q(\chi))=0$, for all $i\le r$, and 
all $q\in \PP(\chi)$ with $q\ne p$. 
\end{claim}

The proof is by induction on $r$, with the case $r=0$ clear;  
indeed, $\sV_q(\chi)\ne \emptyset$ implies 
$\beta_0(\SR,\sV_q(\chi))=0$. 

Assume the claim is true for $r-1$.  
Suppose condition $(\dagger)_{r}$ is satisfied, 
and $T^{\le r}_a=0$, for all $a\in\k\setminus\set{0,1}$. 
Using \eqref{eq:mono eqchain}, we deduce that  
\begin{equation}
\label{eq:hi tl} \tag*{$(*)_r$}
H_i(T_L,\Lambda_\chi)=(\Lambda/(t-1))^{e_i}, \quad 
\text{for all $i\le r$}. 
\end{equation}

Fix an element $a$ as above. Since $a\ne 1$, 
the Universal Coefficient Theorem, together with $(*)_r$, 
gives $H_r(T_L,\k_a)=0$. Hence,  by Lemmas \ref{lem:change vars} 
and \ref{lem:aah}, $\beta_r(\SR,\sW_a)=0$. 

Now fix a prime $q\in \PP(\chi)$, $q\ne p$. 
Since $\k=\bar{\k}$, we may find an element 
$a\in \k^{\times}\setminus \set{1}$ with $\ord(a)=q$. 
Clearly, $\sW_a= \sV_q(\chi)$, and so, by the above, 
$\beta_r(\SR,\sV_q(\chi))=0$. 
This finishes the proof of Claim \ref{claim:1}. 

\begin{claim}
\label{claim:2}
Suppose conditions $(\dagger)_r$ and $(\ddagger)_r$ are 
satisfied. Then $T^i_a=0$, for all all $i\le r$, and all 
$a\in\k\setminus\set{0,1}$.
\end{claim}

The proof is by induction on $r$, with the case $r=0$ clear;  
indeed, $H_0(T_L,\k\Z_\chi)=\k$, and so $T^0_a=0$. 

Assume the claim is true for $r-1$, and suppose  
conditions $(\dagger)_{r}$ and $(\ddagger)_{r}$ are satisfied. 
Using \eqref{eq:mono eqchain}, we deduce that $(*)_{r-1}$ holds, 
and  
\begin{equation}
\label{eq:hr tl}
H_r(T_L,\Lambda_\chi)=(\Lambda/(t-1))^{e_r} \oplus 
\bigoplus_{b\in \k\setminus \set{0,1}} T^r_b. 
\end{equation}

Now fix an element $a$ in $\k\setminus \set{0,1}$. 
We need to show $T^r_a=0$. Write 
\begin{equation}
\label{eq:tra}
T^r_a= \bigoplus_{j\ge 1} \big( \L/(t-a)^j\L \big)^{\varepsilon_j}.
\end{equation}
The Universal Coefficient Theorem, together with $(*)_{r-1}$ 
and \eqref{eq:hr tl} yields $H_r(T_L,\k_a)= \bigoplus_{j\ge 1} 
\k^{\varepsilon_j}$. So it is enough to show $H_r(T_L,\k_a)=0$. 
In view of Lemma \ref{lem:change vars}, 
it remains to show that 
\begin{equation}
\label{eq:betar}
\beta_r(\SR, \sW_a)=0.
\end{equation}

By Lemma \ref{lem:supp}, there are two possibilities to 
consider.
\begin{romenum}
\item $\sW_a=\sV_0(\chi)$.  Recall $\sV_0(\chi)\supset \sV_p(\chi)$. 
By Lemma \ref{lem:semi} and condition $(\dagger)_r$, we 
have $\beta_r(\SR, \sW_a)\le \beta_r(\SR, \sV_p(\chi))=0$.  
\item $\sW_a\supset \sV_q(\chi)$, for some prime $q\ne p$. 
Here, there are two sub-cases. 
\begin{enumerate}
\item $q\in \PP(\chi)$. Then 
$\beta_r(\SR, \sW_a)\le \beta_r(\SR, \sV_q(\chi))=0$.  
\item $q\not\in \PP(\chi)$.  Then $\sW_a\supset \sV_q(\chi)
=\sV_0(\chi)\supset \sV_p(\chi)$, and so 
$\beta_r(\SR, \sW_a)\le \beta_r(\SR, \sV_p(\chi))=0$.  
\end{enumerate}
\end{romenum}

Thus, in all cases equality \eqref{eq:betar} holds.  
This finishes the proof of Claim \ref{claim:2}, and 
thereby ends the proof of the Theorem. 
\end{proof}

\begin{remark}
\label{cor:triv res}
Let $X$ be a connected CW-complex of finite type. Let 
$\nu\colon \pi_1(X)\surj A$ be an epimorphism onto an  
abelian group $A$. Denote by 
$X^{\nu}$ the Galois $A$-cover of $X$ corresponding 
to $\ker (\nu)$. For a coefficient field $\k$, consider 
the following properties of $X^{\nu}$:
\begin{enumerate}
\item \label{f1}
$H_{\le r}(X^{\nu}, \k)$ has trivial $A$-action.
\item \label{f2}
$\dim_{\k} H_{\le r}(X^{\nu}, \k)<\infty$.
\item \label{f3}
$\nu_{\k}\not\in \bigcup_{i=1}^r \RR_1^i (X, \k)$ (when $A=\Z$).
\end{enumerate}

It is readily seen that \eqref{f1} $\Rightarrow$ \eqref{f2}. 
For $A=\Z$, the implications \eqref{f1} $\Rightarrow$ 
\eqref{f3} $\Rightarrow$ \eqref{f2} follow from 
\cite[Proposition 9.4]{PS-eqchain}.   As illustrated by the 
following example (inspired by \cite[Example 8.5]{PS-eqchain}), 
neither of these implications can be reversed. 
\end{remark}

\begin{example}
\label{ex:121}
Consider the cover $T_{\G}^{\chi}\to T_{\G}$ defined 
by the weighted graph
\begin{equation}
\label{eq:path}
\xygraph{
[]!{(0,0)}
*D(3){1}{\cir<2pt>{}}
-@[]!{(1,0)}
*D(3){2}{\cir<2pt>{}}
-@[]!{(2,0)}
*D(3){1}{\cir<2pt>{}}
}
\end{equation}
Clearly, $T_{\G}\simeq S^1\times (S^1\vee S^1)$. If 
$\ch \k \ne 2$, then $\chik\not\in \RR_1^1 (T_{\G}, \k)$, 
by Corollary \ref{thm:res artin}, yet the $\Z$-action on 
$H_1(T_{\G}^{\chi}, \k)$ is non-trivial, by Theorem \ref{thm:mono test}. 
If $\ch \k = 2$, then 
$\dim_{\k} H_{\le 1}(T_{\G}^{\chi}, \k)<\infty$, by 
Corollary \ref{cor:hom res}, but $\chik\in \RR_1^1 (T_{\G}, \k)$, 
again by Corollary \ref{thm:res artin}.  Note also that the 
corresponding Artin kernel, $N_{\chi}=\pi_1(T_{\G}^{\chi})$, 
is finitely presented, as may be checked by using  
test \eqref{test2}  from \S\ref{subsec:zcovers}.
\end{example}

\section{Cohomology ring and finiteness properties}
\label{sect:coho}

We are now in position to compute the (truncated) 
cohomology ring of an arbitrary Galois $\Z$-cover 
of a toric complex, in the case when the monodromy 
action is trivial (up to a fixed degree).

\subsection{Cohomology ring of a toric cover}
\label{subsec:coho ak}

Let $T_L$ be a toric complex, $\chi\colon G_L \surj \Z$ an 
epimorphism, and $\pi\colon T^{\chi}_L \to T_L$ the corresponding 
Galois $\Z$-cover. As before, we denote by $\chi_\k\in H^1(T_L,\k)$ 
the cohomology class determined by $\chi$.  Let $(\chi_\k)$ be the 
ideal of $H^*(T_L,\k)$ it generates. By Proposition 
\ref{prop:coho zcover}\eqref{coho1}, the induced homomorphism 
between cohomology rings, $\pi^*\colon H^*(T_L,\k) \to 
H^*(T^{\chi}_L,\k)$, factors through a ring map,
\begin{equation}
\label{eq:pi map}
\bar\pi^*\colon H^*(T_L,\k)/(\chi_\k) \to H^*(T^{\chi}_L,\k). 
\end{equation}

\begin{theorem}
\label{thm:ak coho}
Let $\k$ be a field, with $\ch \k=p$. Fix an integer $r\ge 1$, 
and suppose 
\begin{enumerate}
\item[$(\dagger)$] 
$\beta_i(\SR,\sV_p(\chi))=0$,
\item[$(\ddagger)$] 
$\beta_i(\SR,\sV_q(\chi))=0$, for all $q\in \PP(\chi)$ 
with $q\ne p$, 
\end{enumerate}
for all $i\le r$.  Then  $\bar\pi^*\colon 
H^{\le r}(T_L,\k)/(\chi_\k) \to H^{\le r}(T^{\chi}_L,\k)$ 
is a ring isomorphism. 
\end{theorem}

\begin{proof}
By Theorem \ref{thm:mono test}, the hypothesis is equivalent 
to $H_{\le r} (T_L,\k\Z_{\chi})$ having trivial $\Z$-action.  
The conclusion follows from Proposition 
\ref{prop:coho zcover}\eqref{coho2}.
\end{proof}

In the particular case of Bestvina-Brady covers, 
$\V_q(\nu)=\V_0(\nu)=\V$, for all primes $q$. In view of 
Proposition \ref{prop:hochster aah}, we recover from 
Theorem \ref{thm:ak coho} the following result of Leary 
and Saadeto\u{g}lu (see Theorem 13 from \cite{LS}). 

\begin{corollary}
\label{cor:bb ls}
Let $\pi\colon T_L^{\nu}\to T_L$ be the Bestvina-Brady cover 
associated to $L$. If $\widetilde{H}_{<r}(L,\k)=0$, 
then $\bar\pi^*\colon H^{\le r}(T_L,\k)/(\nu_\k) 
\to H^{\le r}(T^{\nu}_L,\k)$ is a ring isomorphism.
\end{corollary}

The case when $L$ is a simply-connected flag complex, $r=2$, 
and $\k=\Q$ was first proved in \cite[Theorem 1.3]{PS-bb}, 
by completely different methods. 

\subsection{Finiteness properties of Artin kernels}
\label{subsec:coho bb}

Recall that a group $G$ is of type $\FP_r$ 
($r\le \infty$) if there is a projective $\Z{G}$-resolution 
$P_{\bullet} \to \Z$ of the trivial $G$-module $\Z$, with 
$P_i$ finitely generated for all $i\le r$. The $\FP_r$ 
condition obviously implies that the homology groups 
$H_i(G, \Z)$ are finitely generated, for all $i\le r$, but 
the converse is far from true, in general. 

When coupled with the main result from \cite{MMV}, as restated 
in a more convenient form in \cite{BG}, our Corollary \ref{cor:hom res} 
yields a remarkable property of Artin kernels:  the fact that, within 
this class of groups, the finiteness property $\FP_r$ may
be detected by the corresponding finiteness property for 
homology with trivial field coefficients.

\begin{theorem}
\label{thm=fptriv}
An Artin kernel, $N_{\chi}=\ker (\chi\colon G_{\G} \surj \Z)$, is 
of type $\FP_r$ if and only if $\dim_{\k} H_{\le r} (N_{\chi}, \k)< \infty$,
for any field $\k$.
\end{theorem}

\begin{proof}
Set $L=\Delta_{\G}$ and $\sW=\V_0(\chi)$.  To prove the 
non-trivial implication, assume 
$\dim_{\k} H_{\le r} (N_{\chi}, \k)< \infty$, for any field $\k$. 
Then, by  Corollary \ref{cor:hom res}\eqref{bns3}, 
\[
\widetilde{H}_{i}(\lk_{L_{\sW}} (\sigma), \k)=0,
\]
for all $\sigma\in L_{\V \setminus \sW}$ and $i \le r-1-\abs{\sigma}$.  
Consequently, $\widetilde{H}_{i}(\lk_{L_{\sW}} (\sigma), \Z)=0$,
for all $\sigma$ and $i$ as above. By \cite[Theorem 14]{BG}, 
the group $N_{\chi}$ is of type $\FP_r$.
\end{proof}

For Bestvina-Brady covers and groups, one can say more. 

\begin{corollary}
\label{cor:bb mono}
Let $L$ be a finite simplicial complex on vertex set $\sV$, and 
let $\nu\colon G_L\to \Z$ be the homomorphism sending 
each generator $v\in \sV$ to $1$. For a field $\k$, and 
an integer $r\ge 1$, the following are equivalent:
\begin{romenum}
\item \label{bb1}
The $\k\Z$-module $H_{i} (T_L^{\nu},\k)$ is trivial, 
for all $i\le r$.
\item \label{bb2}
The $\k$-vector space $H_{i} (T_L^{\nu},\k)$ is 
finite-dimensional, for all $i\le r$.
\item \label{bb3}
$\nu_{\k} \not\in \RR_1^i (T_L, \k)$, for all $i\le r$.
\item \label{bb4}
$\widetilde{H}_{i}(L,\k)=0$, for all $i< r$. 
\end{romenum}
If, in addition, $L=\Delta_\Gamma$ is a flag complex, 
conditions \eqref{bb1}--\eqref{bb4} hold over fields $\k$ 
of arbitrary characteristic if and only if $N_{\Gamma}=\ker(\nu)$ 
is of type $\FP_r$. 
\end{corollary}

\begin{proof}
For the implications \eqref{bb1} $\Rightarrow$ \eqref{bb3} 
$\Rightarrow$ \eqref{bb2}, see Remark \ref{cor:triv res}. 
The implication \eqref{bb2} $\Rightarrow$ \eqref{bb1} follows 
from Corollary \ref{cor=bbcover}, while the equivalence 
\eqref{bb2} $\Leftrightarrow$ \eqref{bb4} follows from 
Corollary \ref{cor:hom res}, since $\V_0(\nu)=\V$. 
Finally, the claim about flag complexes follows from 
Theorem \ref{thm=fptriv}.
\end{proof}

\section{Holonomy Lie algebra}
\label{sect:holo}

In this section, we study a certain graded Lie algebra, 
$\h(A)$, attached to a strongly graded-commutative algebra $A$. 
In the process, we relate the non-resonance properties 
of an element $a\in A^1$ to the graded ranks of 
the metabelian quotient of $\h(A/aA)$. 

\subsection{Holonomy Lie algebra of an algebra}
\label{subsec:holo lie}
Let $A$ be a connected, graded, graded-commutative algebra 
over a field $\k$, with graded pieces $A^i$, $i\ge 0$. We shall 
assume that $\dim_{\k} A^1 < \infty$, and 
$a^2=0$, for all $a\in A^1$.  This last condition (which is 
automatically satisfied if $\ch\k\ne 2$) insures that the 
multiplication map in degree $1$ descends to a linear map 
$\mu\colon A^1 \wedge A^1 \to A^{2}$.  

Let $A_i=(A^i)^{\#}$ be the dual $\k$-vector space, and 
let $\Lie(A_1)$ be the free Lie algebra on $A_1$, graded 
by bracket length.  Let $\nabla \colon A_2 \to 
A_1\wedge A_1 = \Lie_2(A_1)$ be the dual of $\mu$. 
In the spirit of K.-T.~Chen's approach from \cite{Ch77}, 
define the {\em holonomy Lie algebra} of $A$ as the quotient 
of the free Lie algebra on $A_1$ by the ideal generated 
by the image of the comultiplication map,
\begin{equation}
\label{eq:holo lie def}
\h(A)=\Lie(A_1)/(\im \nabla). 
\end{equation}

Clearly, $\h=\h(A)$ inherits a natural grading from the 
free Lie algebra, compatible with the Lie bracket; 
let $\h_s$ be the $s$-th graded piece.  By construction, 
$\h$ is a finitely presented graded Lie algebra, with 
generators in degree $1$, and relations in degree $2$.

Note that $\h(A)$ depends only on the degree $2$ truncation 
of $A$,  namely $A^{\le 2}=\bigoplus_{i\le 2} A^i$.   For the 
purpose of defining the holonomy Lie algebra of $A$, we may 
assume $A$ is generated in degree $1$.  Indeed, if we let 
$E=\bigwedge A^1$ be the exterior algebra on $A^1$, and 
we set $\bar{A}=E/(K+ E^{\ge 3})$, where 
$K=\ker(\mu)$, then the algebra $\bar{A}$ is generated by 
$\bar{A}^1=A^1$, and clearly $\h(\bar{A})=\h(A)$.

So assume $A$ is generated in degree $1$, i.e, the morphism 
$q_A\colon E=\bigwedge A^1 \to A$ extending the 
identity on $A^1$ is surjective.  In this case, the 
map $\nabla$ is injective. Thus, we may view 
$A_2$ as a subspace of $A_1\wedge A_1$, 
and arrive at the following identifications:
\begin{equation}
\label{eq:holo123}
\h_1=A_1, \quad  \h_2 = (A_1 \wedge A_1)/A_2, 
\quad  \h_3=\Lie_3(A_1)/[A_1,A_2].
\end{equation}
\subsection{Homological algebra interpretation}
\label{subsec:tor}
Since $\h$ is generated in degree $1$, 
its derived Lie subalgebra, $\h'$, coincides 
with $\h_{\ge 2}=\bigoplus_{s\ge 2} \h_s$.  Let $\h''$ be 
the second derived Lie subalgebra of $\h$, and let 
\begin{equation}
\label{eq:g lie def}
\g(A)=\h(A)/\h''(A)
\end{equation}
be the maximal metabelian quotient Lie algebra of $\h$.  
It is readily checked that $\g'=\h'/\h''$. Moreover, 
$\g_s=\h_s$, for $s\le 3$. 

Viewing $A$ as an $E$-module via $q_A$, and 
$\k$ as the trivial $E$-module $E/E^{\ge 1}$, we may form 
the bigraded vector space $\Tor^E_{*}(A,\k)^{\phantom{E}}_{*}$, 
where the first grading comes from a free resolution of $A$ over 
$E$, and the second one comes from the grading on $E$. Then:
\begin{equation}
\label{eq:fl}
\g_s(A)^{\sharp} =  \Tor^E_{s-1}(A,\k)_{s},
\end{equation}
for all $s\ge 2$.   As noted in Proposition 2.3 from 
\cite{PS-artin}, equality \eqref{eq:fl} follows easily from work 
of Fr\"oberg and L\"ofwall (Theorem 4.1(ii) in \cite{FL}). 

\subsection{Resonance and holonomy}
\label{subsec:res holo}
As above, let $A$ be an algebra with $A^1$ finite-dimen\-sional, 
$a^2=0$ for all $a\in A^1$, and $A$ generated in degree $1$.  
Let $a$ be a non-zero element in $A^1$, and consider the 
quotient algebra $B=A/aA$.  The next result equates the 
graded ranks of $\g'(A)$ and $\g'(B)$ in a certain range, 
prescribed by the non-resonance properties of $a$. 

\begin{theorem}
\label{thm:hab}
Suppose $a\notin \bigcup_{i=1}^r \RR^i_1(A)$, 
for some $r\ge 1$. Then 
\begin{equation*}
\label{eq:hab}
\g_{s}(A) \cong \g_{s}(B), \quad \text{for $2\le s\le r+1$}.
\end{equation*}

\end{theorem}

\begin{proof}
Set $F=\bigwedge B^1$. 
Let  $\phi\colon A\to B$ be the projection map, 
and let $\psi\colon E\to F$ be the extension of 
$\phi^1\colon A^1\to B^1$ to exterior algebras.  
We then have a commuting square  
\begin{equation}
\label{eq:sq}
\xymatrix{
E \ar@{>>}^{\psi}[r] \ar@{>>}^{q_A}[d] 
& F\ar@{>>}^{q_B}[d] \\
A \ar@{>>}^{\phi}[r] & B
}
\end{equation}
with $F=E/aE$ and $A\otimes_E F=B$. This square 
gives rise to a ``change of rings" spectral sequence, 
\begin{equation}
\label{eq:ss}
E^2_{s,t}=\Tor^F_s (\Tor^E_t (A,F) ,\k) \Rightarrow 
\Tor^E_{s+t} (A,\k).
\end{equation}
Note that the $E^2$-term has an extra grading, $_{*}E^2_{s,t}$, 
coming from the degree grading on $E$.  This extra 
grading is preserved by the differentials, and is compatible 
with the internal grading on $\Tor^E(A,\k)$. 

A free resolution of the $E$-module $F$ is given by 
\begin{equation}
\label{eq:f res}
\xymatrix{
\cdots \ar[r] & E[t] \ar^(.4){a}[r] & E[t-1] \ar[r] & \cdots \ar[r]&
E[1] \ar^(.55){a}[r] & E \ar^{\psi}[r] \ar[r] & F \ar[r] & 0}.
\end{equation}
Here, the free module in position $t$ is regraded as $E[t]$, 
with $E[t]^i=E^{i-t}$; this is done so that all boundary maps 
have degree $0$. Therefore, $\Tor^E_t (A,F)$ is the 
homology in degree $t$ of the chain complex 
\begin{equation}
\label{eq:ef res}
\xymatrix{ \cdots \ar[r] &A[t] \ar^(.4){a}[r]& A[t-1] \ar[r] & \cdots 
\ar[r]& A[1] \ar^(.55){a}[r] & A \ar^{\phi}[r]& B \ar[r] & 0},
\end{equation} 
In internal  degree $i$, this homology group 
coincides with $H^{i-t}(A,\cdot a)$, for $t>0$. Thus,
\begin{equation}
\label{eq:tor h}
\Tor^E_t (A,F)_i = H^{i-t}(A,\cdot a), \quad \text{for $t>0$} .
\end{equation}

By assumption, $a\notin \RR^{i}_1(A)$, that is, $H^{i}(A,\cdot a)=0$, 
for $i\le r$.  Hence, by \eqref{eq:tor h}, we have 
$\Tor^E_t (A,F)_{i}=0$ for $i\le r+t$ and $t>0$. It follows that 
$_{i}E^2_{s,t}=0$, for all $s\ge 0$, $t>0$ and $i\le r+1$, 
and so 
\begin{equation}
\label{eq:tor ab}
\Tor^E_{s-1}(A,\k)_{s} = _{s\!}E^{\infty}_{s-1,0}
=_{s\!}E^2_{s-1,0}= \Tor^F_{s-1}(B,\k)_{s},
\end{equation}
for all $1\le s\le r+1$. Invoking \eqref{eq:fl} finishes the proof.
\end{proof}

\begin{corollary}
\label{cor:res holo}
Let $a\in A^1$, and set $B=A/aA$.  
\begin{enumerate}
\item \label{c1}  If $a\notin \RR^1_1(A)$, then $\h_2(A)=\h_2(B)$.
\item \label{c2}  If $a\notin \RR^1_1(A)\cup \RR^2_1(A)$, then 
$\h_3(A)=\h_3(B)$. 
\end{enumerate}
\end{corollary}

Let $\LL_{m}=\Lie(\k^m)$ be the free Lie $\k$-algebra on $m$ generators, 
graded by bracket length.

\begin{corollary}
\label{cor:res chen}
Suppose  $A^1= \k^n$, $n\ge 1$, and $A^{\ge 3}=0$.  
If $\RR^1_1(A)\cup \RR^2_1(A)\ne A^1$, then 
$\g' (A)\cong \LL'_{n-1}/\LL''_{n-1}$.
\end{corollary}

\begin{proof}
Pick $a\in A^1 \setminus (\RR^1_1(A)\cup \RR^2_1(A))$, 
and set $B=A/a A$. Non-resonance of $a$ in degree $1$ 
implies $a\ne 0$; hence, $B^1=\k^{n-1}$. Non-resonance 
of $a$ in degree $2$ implies $A^2= a A^1$; in particular, $A$ is 
generated in degree $1$, and $B^{\ge 2}=0$. Hence, $\h (B)=\LL_{n-1}$, 
by definition \eqref{eq:holo lie def}.  The claim then follows from 
Theorem \ref{thm:hab}.
\end{proof}

\section{Lie algebras associated to Artin kernels}
\label{sect:lie ak}

In this section, we study three graded Lie algebras associated 
to an Artin kernel $N_{\chi}=\ker(\chi\colon G_L\surj \Z)$: the 
associated graded Lie algebra $\gr(N_{\chi})$, the Chen Lie algebra 
$\gr(N_{\chi}/N''_{\chi})$, and the holonomy Lie algebra $\h(N_{\chi})$. 

\subsection{Associated graded and Chen Lie algebras}
\label{subsect:gr}

We start with a classical construction from the 1930s, 
due to P.~Hall and W.~Magnus, see \cite{MKS}.  
Let $G$ be a group. The lower central series (LCS) 
of $G$  is defined inductively by $\gamma_1 G=G$ 
and $\gamma_{k+1}G =(\gamma_k G,G)$, where 
$(x,y)=xyx^{-1}y^{-1}$. The {\em associated 
graded Lie algebra}, $\gr(G)$, is the direct sum 
of the successive LCS quotients, 
\begin{equation}
\label{eq:gr g}
\gr(G)=\bigoplus\nolimits_{k\ge 1} 
\gamma_k G/ \gamma_{k+1} G,
\end{equation} 
with Lie bracket induced from the group commutator.
For a field $\k$, 
we shall write $\gr(G)\otimes \k= 
\bigoplus_{k\ge 1} \gr_k(G)\otimes \k$. 

By construction, the Lie algebra $\gr(G)$ is generated 
by $\gr_1(G)=G_{\ab}$. Thus, if $G$ is a finitely generated 
group, then $\gr(G)$ is a finitely generated Lie algebra. 
Moreover, the derived Lie subalgebra, $\gr'(G)$, coincides 
with $\bigoplus\nolimits_{k\ge 2} \gr_k(G)$. 

In \cite{Ch51}, K.-T.~Chen introduced a useful variation 
on this theme. Let $G'=\gamma_2 G$ be the derived group, 
and $G''=(G')'$ the second derived group.  Note that 
$H_1(G,\Z)=G/G'$ is the maximal abelian quotient of $G$, 
whereas $G/G''$ is the maximal metabelian quotient.  
The {\em Chen Lie algebra}\/ of $G$ is simply $\gr(G/G'')$.  
Though a coarser invariant that $\gr(G)$, the Chen Lie algebra 
captures some subtle phenomena, in its own distinctive 
manner, more closely tied to commutative algebra. 
We refer to \cite{PS-chen} for more on this subject. 

\subsection{Associated graded of an Artin kernel}
\label{subsect:gr ak}

Let $\G$ be a finite simple graph, with vertex set $\V$, 
and denote by $G=G_{\G}$ the corresponding right-angled 
Artin group. Let $\chi\colon G\surj \Z$ be an epimorphism, 
and set $N=\ker(\chi)$. Denoting by $\iota\colon N\to G$ 
the inclusion map, we have a split exact sequence
\begin{equation}
\label{eq:iota nu}
\xymatrix{1\ar[r] &  N \ar^{\iota}[r] &
G \ar^{\chi}[r] & \Z \ar[r] & 1}.
\end{equation}

\begin{lemma}
\label{lem:fg}
Suppose $H_1(N,\Q)$ has trivial $\Q\Z$-action. Then:
\begin{enumerate}
\item \label{nq1}  $N$ is finitely generated.  
\item \label{nq2}  $H_1(N,\k)$ has trivial $\k\Z$-action, 
for any field $\k$.
\item  \label{nq3}  $H_1(N,\Z)$ has trivial $\Z\Z$-action.
\item  \label{nq4}  The restriction of $\iota$ to derived subgroups, 
$\iota'\colon N'\to G'$, is an isomorphism. 
\item  \label{nq5}  The induced homomorphism $\iota_*\colon H_1(N,\Z) \to 
H_1(G,\Z)$ is injective.  In particular, $H_1(N,\Z)$ is a free abelian group, 
of rank $\abs{\sV}-1$. 
\end{enumerate}
\end{lemma}

\begin{proof}
\eqref{nq1}
By Theorem \ref{thm:mono test} and Corollary \ref{cor:hom res}, 
the induced subgraph $\Gamma_{\V_0(\chi)}$ is connected and 
dominant.  Hence, by \cite{MV}, the group $N$ is finitely generated.  

\eqref{nq2}
By Theorem \ref{thm:mono test}, our hypothesis means that 
$\beta_i(\Q \langle L\rangle, \sW)=0$, for $i=0, 1$, and 
\begin{equation}
\label{eq:beta0}
\text{$\sW=\sV_0(\chi)$, or $\sW=\sV_q(\chi)$, with $q\in \PP(\chi)$}.
\end{equation}

Now let $\k$ be a field of characteristic $p$. Again 
by Theorem \ref{thm:mono test}, we have to check that 
$\beta_i(\SR, \sW)=0$, for $i=0,1$, where either $\sW=\sV_p(\chi)$, 
or $\sW=\sV_q(\chi)$, with $q\in \PP(\chi)$ and $q\ne p$. 
Using Remark \ref{rem:beta} and formula \eqref{eq:beta1}, 
we see that it is enough to check $\beta_i(\Q \langle L\rangle, \sW)=0$, 
for $i=0,1$, and $\sW$ as just above.  There are $3$ cases 
to consider:
\begin{itemize}
\item $\sW=\sV_q(\chi)$, with $q\in \PP(\chi)$, $q\ne p$.
\item $\sW=\sV_p(\chi)$, and $p\in \PP(\chi)$.
\item $\sW=\sV_p(\chi)$, and $p\notin \PP(\chi)$, 
in which case $\sV_p(\chi)=\sV_0(\chi)$, by 
definition \eqref{eq:pnu}.
\end{itemize}
Clearly, all three cases are covered by \eqref{eq:beta0},   
and we are done.

\eqref{nq3} Follows from \eqref{nq1}, \eqref{nq2}, and  
the Universal Coefficient Theorem.

\eqref{nq4}
Plainly, $\iota (N')\subseteq G'$.  The triviality of the $\Z\Z$-action 
on $N/N'=H_1(N,\Z)$ forces $\iota (N')=G'$. Thus, 
$\iota'\colon N'\to G'$ is an isomorphism.  

\eqref{nq5}
It follows from \eqref{nq3} and \cite[Lemma 3.4]{FR}
that \eqref{eq:iota nu} remains exact upon abelianization.  
Recalling that $H_1(G,\Z)$ is a free abelian 
group of rank $\abs{\sV}$ finishes the proof. 
\end{proof}

For a homomorphism $\alpha\colon G\to H$, let 
$\bar\alpha\colon G/G'' \to H/H''$ be the induced 
homomorphism on maximal metabelian quotients, 
and $\gr(\alpha)\colon \gr(G)\to \gr(H)$ the induced 
morphism of graded Lie algebras. The proof of the 
next result is similar to (and generalizes) the proofs 
of Propositions 4.2 and 5.4, and Theorem 5.6 from 
\cite{PS-bb}.  

\begin{prop}
\label{prop:gr split}
Suppose $H_1(N,\Q)$ has trivial $\Q\Z$-action.
Then, we have split exact sequences of graded Lie algebras,    
\begin{align}
\label{eq:fr}
&\qquad \xymatrixcolsep{14pt}
\xymatrix{0\ar[r] &  \gr(N)\ar^{\gr(\iota)}[rr] &&
\gr(G)  \ar^{\gr(\chi)}[rr] && \gr(\Z)  \ar[r] & 0},
\\
\label{eq:fr chen}
&\xymatrixcolsep{14pt}
\xymatrix{0\ar[r] &  \gr(N/N'')\ar^{\gr(\bar\iota)}[rr] &&
\gr(G/G'') \ar^(.55){\gr(\bar\chi)}[rr] && \gr(\Z)  \ar[r] & 0}.
\end{align}
\end{prop}

\begin{proof} 
By Lemma \ref{lem:fg}\eqref{nq3}, the abelianization 
$N/N'=H_1(N,\Z)$ is a trivial $\Z\Z$-module. 
Applying the $\gr$ functor to \eqref{eq:iota nu}, 
and making use again of Lemma 3.4 from Falk 
and Randell \cite{FR}, yields \eqref{eq:fr}. 

By Lemma \ref{lem:fg}\eqref{nq4}, $\iota (N')=G'$.  
It follows that $\iota (N'')=G''$.  Hence, there is  
an exact sequence $1\to  N/N'' \xrightarrow{\,\bar\iota\,}
G/G'' \xrightarrow{\,\bar\chi\,} \Z \to 1$.  Applying the $\gr$ 
functor as above yields \eqref{eq:fr chen}. 
\end{proof}

Consequently, $\gr'(N)\cong \gr'(G)$ and $\gr'(N/N'')\cong \gr'(G/G'')$, 
as graded Lie algebras. Using the computation of $\gr(G)$ and 
$\gr(G/G'')$ from \cite{PS-artin}, we obtain the following immediate 
corollary, which generalizes Theorems 5.1 and 5.2 from \cite{PS-bb}. 

\begin{corollary}
\label{cor:gr ak}
Suppose $H_1(N,\Q)$ has trivial $\Z$-action.
Then, both $\gr(N)$ and $\gr(N/N'')$ are torsion-free, with 
graded ranks, $\phi_k$ and $\theta_k$,  given by 
\begin{equation}
\label{eq:lcs ranks artin}
\prod_{k=1}^{\infty}(1-t^k)^{\phi_k}=\frac{P(-t)}{1-t},
\end{equation}
where $P(t)=\sum_{k\ge 0} f_k(\G) t^k$ is the {\em clique 
polynomial} of $\G$, with $f_k(\G)$ equal to the number 
of $k$-cliques of $\G$, and 
\begin{equation}
\label{eq:chen ranks artin}
\sum_{k=2}^{\infty} \theta_k t^{k} = 
Q \Big(\frac{t}{1-t}\Big),
\end{equation}
where $Q(t)=\sum_{j\ge 2} 
\left(\sum_{\sW\subset \sV\colon\,  \abs{\sW}=j } 
\tilde{b}_0(\G_{\sW})\right) t^j$ 
is the {\em cut polynomial} of $\G$.
\end{corollary}

\subsection{Holonomy Lie algebra of an Artin kernel}
\label{subsect:holo ak}

Suppose $G$ is a finitely generated group.  Then 
$A=H^*(G;\k)$ is a connected, graded, graded-commutative 
$\k$-algebra, with $A^1$ finitely generated.   Assuming 
$\ch\k\ne 2$, we may define the holonomy Lie algebra 
of $G$, with coefficients in $\k$, to be $\h(G,\k)= \h(A)$; 
we will simply write  $\h(G)=\h(G,\Q)$. 

The obvious identification $A_1=\gr_1(G)\otimes \Q$ 
extends to a Lie algebra map, $\Lie(A_1) \surj \gr(G)\otimes \Q$, 
which in turn factors through an epimorphism 
$\Psi_G \colon  \h(G) \surj \gr (G)\otimes \Q $; 
see \cite{PS-chen} for further details and references. 

Consider now a right-angled Artin group, 
$G=G_{\Gamma}$.  In this case, it is known that  
$\Psi_G \colon  \h(G) \to \gr (G)\otimes \Q $ is 
an isomorphism; see \cite{PS-artin}. 
Let $\chi\colon G\surj \Z$ be an epimorphism, 
and $N= \ker (\chi)$ the corresponding Artin kernel. 

\begin{lemma}
\label{lem:holo surj}
Suppose $H_1(N,\Q)$ has trivial $\Q\Z$-action.
Then:
\begin{enumerate}
\item \label{iota1} $\h_1(\iota)\colon \h_1(N) \to \h_1(G)$ is injective. 
\item \label{iota2} $\h_k(\iota)\colon \h_k(N) \to \h_k(G)$ is surjective, 
for all $k\ge 2$.
\end{enumerate}
\end{lemma}

\begin{proof}
Part \eqref{iota1}.  The map $\h_1(\iota)$ may be identified 
with the induced homomorphism 
$\iota_*\colon H_1(N,\Q)\to H_1(G,\Q)$, which is injective, 
by triviality of the $\Q\Z$-action on $H_1(N,\Q)$. 

Part \eqref{iota2}. Consider the commuting diagram 
\begin{equation}
\label{eq:gr hol cd}
\xymatrix{\h(N)\ar^{\h(\iota)}[r] \ar@{>>}^{\Psi_N} [d] 
& \h(G)\ar@{>>}^{\Psi_G}[d] \\ 
\gr(N) \otimes \Q\ar^{\gr(\iota)}[r] & \gr(G) \otimes \Q
}
\end{equation}
and fix a degree $k\ge 2$. By Proposition \ref{prop:gr split}, 
$\gr_k(\iota)$ is an isomorphism.  Thus, $\h_k(\iota)$ is surjective. 
\end{proof}

The next Theorem identifies (under certain conditions) 
the holonomy Lie algebra $\h(N)$ as a Lie subalgebra 
of $\h(G)$.  

\begin{theorem}
\label{thm:holo split}
Let $N=\ker(\chi\colon G\surj \Z)$. Suppose
$H_i(N,\Q)$ has trivial $\Q\Z$-action, for $i=1,2$.  
Then, we have a split exact sequence of graded Lie algebras,    
\begin{equation}
\label{eq:holo ex}
\xymatrix{0\ar[r] &  \h(N)\ar^{\h(\iota)}[r] &
\h(G)  \ar^{\h(\chi)}[r] & \h(\Z)  \ar[r] & 0}.
\end{equation}
In particular, the restriction of $\h(\iota)$ to derived 
Lie subalgebras, $\h'(\iota)\colon \h'(N)\to \h'(G)$, is an 
isomorphism of graded Lie algebras. 
\end{theorem}

\begin{proof}
Let $A=H^{\le 2}(G,\Q)$, so that $\h(G)=\h(A)$ . 
Consider the element $a=\chi_{_{\Q}}\in A^1$, 
and set $B=A/aA$.  By Theorem \ref{thm:ak coho}, 
$B\cong H^{\le 2}(N,\Q)$.  Therefore, $\h(N)=\h(B)$.  

Since $H_{\le 2}(N,\Q)$ has trivial $\Q\Z$-action, 
$a$ does not belong to $\RR^1_1(G, \Q)\cup \RR^2_1(G, \Q)$;
see Remark \ref{cor:triv res}.
Hence, by Corollary \ref{cor:res holo}, 
$\h_2(G)=\h_2(N)$ and $\h_3(G)=\h_3(N)$. 
It follows from Lemma \ref{lem:holo surj} that 
$\h_2(\iota)$ and $\h_3(\iota)$ are isomorphisms. 

Set $n=b_1(G)$, and pick bases $\{x_1,\dots ,x_{n-1}\}$ 
for $B_1$ and $\{x_1,\dots ,x_{n-1}, y\}$ for $A_1$.  
Since $\h_2(\iota)\colon (B_1 \wedge B_1)/B_2 \to 
(A_1 \wedge A_1)/A_2$ is an isomorphism, we may 
identify $A_2/B_2$ with $A_1 \wedge A_1
/B_1 \wedge B_1$, which is a vector space with basis 
the cosets represented by 
$\{ x_1 \wedge y,\dots ,  x_{n-1}  \wedge y\}$. 
Thus, we may decompose $A_2$ as 
\[
A_2 = B_2 \oplus \spn \{ x_1 \wedge y- z_1,\dots , 
x_{n-1}  \wedge y-z_{n-1}\},
\]
where $z_i\in B_1 \wedge B_1$. 
Now define a linear map, 
\[
\alpha\colon B_1 \to B_1 \wedge B_1, \quad 
\alpha(x_i )= z_i,
\]
and extend it to a degree~$1$ Lie algebra derivation, 
$\tilde\alpha\colon \Lie(B_1) \to \Lie(B_1)$. Given 
an element $s\in B_2\subset \Lie_2(B_1)$, we have  
$\h_3(\iota) (\tilde\alpha(s)) =[\h_2 (\iota)(s),y]$, and 
this vanishes in $\h_3(A)$. 
Since $\h_3(\iota)$ is an isomorphism, $\tilde\alpha(s)=0$ 
in $\h_3(B)$.  In other words, $\tilde\alpha(B_2)\subset [B_1,B_2]$. 
Hence, $\tilde\alpha$ factors through a degree $1$ Lie derivation, 
$\tilde\alpha\colon \h(B) \to \h(B)$.  Since the graded Lie algebra
$\h(\Z)$ is freely generated by $y$, the Lie algebra $\h(A)$ 
splits as  a semidirect product, 
$\h(A)=\h(B)\rtimes_{\tilde\alpha} \h(\Z)$, and we are done. 
\end{proof}

For the Bestvina-Brady groups $N_{\G}=\ker(\nu\colon 
G_{\G}\surj \Z)$, the above results take a particularly 
simple form.  

\begin{corollary}
\label{cor:bb holo}
Let $\G$ be a finite, connected graph.  If $H_1(\Delta_{\G},\Q)=0$, 
then the inclusion map $\iota\colon N_{\G}\inj G_{\G}$ 
induces a group isomorphism, $\iota' \colon N'_{\G}\isom G'_{\G}$,
and an isomorphism of graded Lie algebras, 
$\h'(\iota)\colon \h'(N_{\G})\isom \h'(G_{\G})$.
\end{corollary}

\begin{proof}
By Corollary \ref{cor:bb mono},  
$H_{\le 2}(N_{\G},\Q)$ has trivial $\Q\Z$-action.
The conclusions follow from Lemma \ref{lem:fg}\eqref{nq4}  
and Theorem \ref{thm:holo split}.
\end{proof}

This corollary generalizes \cite[Lemma 6.2]{PS-bb}, proved 
under the more restrictive hypothesis $\pi_1(\Delta_{\G})=0$.  

\section{Formality properties of Artin kernels}
\label{sect:formal}

In this section we give certain conditions guaranteeing 
the $1$-formality of a finitely-generated Artin kernel 
$N_{\chi}$.  These conditions may be satisfied, even 
when $N_{\chi}$ does not admit a finite presentation. 

\subsection{Malcev Lie algebras and $1$-formality} 
\label{subsect:malcev}

In \cite[Appendix A]{Q}, Quillen defines a {\em Malcev 
Lie algebra}\/ to be a rational Lie algebra $E$, endowed 
with a decreasing, complete, $\Q$-vector space filtration,
$E=F_1E\supseteq F_2E\supseteq  \dots$, 
with the property that $[F_sE, F_rE]\subseteq F_{s+r}E$, 
for all $s,r\ge 1$, and such that the associated graded 
Lie algebra, $\gr (E)$, is generated by $\gr_1(E)$. 

An example of a Malcev Lie algebra is $\wh(G)$, the 
completion with respect to bracket length filtration 
of $\h(G)$, the rational holonomy Lie algebra of a 
finitely generated group $G$. Clearly, $\gr(\wh(G))=\h(G)$.  

Also in \cite{Q}, Quillen associates to a group 
$G$ a pronilpotent, rational Lie algebra, $\M(G)$, called 
{\em the}\/ Malcev Lie algebra of $G$.  This functorial 
construction yields a Malcev Lie algebra with the 
crucial property that $\gr(\M(G)) \cong \gr(G)\otimes \Q$, 
as graded Lie algebras.  

Assume now $G$ is finitely generated. Following 
Sullivan \cite{S}, we say $G$ is {\em $1$-formal}\/ 
if $\M(G)\cong \wh(G)$, as filtered Lie algebras.
Equivalently, $\M(G)$ is filtered Lie isomorphic to 
the degree completion of a quadratic Lie algebra. 
This fact is proved in \cite[Lemma 2.9]{DPS-serre} 
for finitely presented groups, but the 
proof given there works as well in this wider generality. 

If $G$ is $1$-formal, the map $\Psi_G\colon \h(G)\to 
\gr(G) \otimes \Q$ is an isomorphism.  Moreover, there 
is a (non-canonical) filtered Lie isomorphism,
$\kappa_G\colon \M(G) \isom \wh(G)$, 
with the property that $\gr_1(\kappa_G)= 
\id_{G_{\ab}\otimes \Q}$; see \cite[Lemma 2.10]{DPS-serre}, 
with the same proviso as above.  

\subsection{$1$-Formality of {A}rtin kernels}
\label{subsect:formal ak}

In \cite{KM}, Kapovich and Millson showed that all 
finitely generated (in particular, right-angled) Artin groups 
are $1$-formal.  In \cite{PS-bb}, we showed that 
finitely presented Bestvina-Brady groups are 
$1$-formal.  We now generalize this last result 
to a wider class of (not necessarily finitely presented) 
Artin kernels. 

As before, let  $G=G_{\Gamma}$ be a right-angled 
Artin group, $\chi\colon G\surj \Z$ an epimorphism, 
and $N=\ker( \chi)$ the corresponding Artin kernel. 

\begin{theorem}
\label{thm:ker formal}
Suppose
$H_i(N,\Q)$ has trivial $\Q\Z$-action, for $i=1,2$.  
Then $N$ is finitely generated and $1$-formal. 
\end{theorem}

\begin{proof}
Finite generation of $N$ is assured by Lemma \ref{lem:fg}\eqref{nq1}. 
As a functor from finitely generated groups to filtered 
Lie algebras, Malcev completion is right-exact.  
Applying the functor $\M$ to the exact sequence 
\eqref{eq:iota nu}, we get an exact sequence 
of filtered Lie algebras, 
\begin{equation}
\label{eq:malcev}
\xymatrix{\M(N) \ar^{\M(\iota)}[r] &
\M(G) \ar^{\M(\chi)}[r] & \M(\Z) \ar[r] & 0}.
\end{equation}

The exactness of \eqref{eq:fr} and the natural isomorphism 
$\gr( \M(\bullet)) \cong \gr(\bullet) \otimes \Q$ imply that 
$\gr(\M(\iota))$ is injective. Since the filtration on 
$\M(N)$ is complete, we conclude that $\M(\iota)$ 
is injective; that is, sequence \eqref{eq:malcev} is 
also exact on the left. 

Let $\h(\chi) \colon \h(G)\surj \h(\Z)$ be the morphism 
induced by $\chi$ at the level of holonomy Lie algebras. 
Passing to completions, we obtain a filtered Lie morphism, 
$\wh(\chi) \colon \wh(G) \surj \wh(\Z)$, whose kernel, 
denoted by $\K$, we equip with the induced filtration. 

Since $G$ and $\Z$ are $1$-formal, we have filtered 
Lie isomorphisms, $\kappa_G \colon \M(G) \isom\wh(G)$ and
$\kappa_{\Z} \colon \M(\Z) \isom \wh(\Z)$, normalized in 
degree $1$, as explained in \S\ref{subsect:malcev}. 
Clearly, $\gr_{>1}(\Z)\otimes \Q=0$, which implies that 
$F_2 \M(\Z)=0$. Using the normalization property, 
we infer that $\kappa_{\Z}\circ \M(\chi)= \wh(\chi)\circ \kappa_G$.
We then have the following commuting diagram in the 
category of filtered Lie algebras:
\begin{equation}
\label{eq:malhol}
\xymatrix{0\ar[r]& \M(N)\ar[d] \ar^{\M(\iota)}[r] &
\M(G)\ar^{\kappa_{G}}[d] \ar^{\M(\chi)}[r] 
& \M(\Z)\ar^{\kappa_{\Z}}[d] \ar[r] & 0\\
0\ar[r]& \K \ar[r] &
\wh(G)  \ar^{\wh(\chi)}[r] & \wh(\Z)  \ar[r] & 0
}
\end{equation}

Since $\gr(\M(\iota))$ is injective, the filtration of $\M(N)$ is 
induced from $\M(G)$. Hence, $\M(N)\cong \K$, as filtered 
Lie algebras. Now note that $\K$ is the kernel of 
$\h(\chi) \colon \h(G) \to \h(\Z)$, completed with 
respect to degree filtration. 
By Theorem \ref{thm:holo split}, 
$\ker (\h(\chi))=\h(N)$. Hence, $\K\cong \wh(N)$, 
as filtered Lie algebras, and we are done. 
\end{proof}

\subsection{Non-finitely presented, $1$-formal groups}
\label{subsec:nonfp formal}

We now apply the above machinery to the Bestvina-Brady 
groups $N_{\Gamma}$.  In \cite[Proposition 6.1]{PS-bb}, 
we proved the following:  If $\Delta_{\Gamma}$ is 
simply-connected (equivalently, if $N_{\Gamma}$ is 
finitely presented), then $N_{\Gamma}$ is $1$-formal. 
We may strengthen that result, as follows.

\begin{corollary}
\label{cor:bb formal}
Let $\Gamma$ be a finite, connected graph.
If $H_1(\Delta_{\Gamma},\Q)=0$, then 
$N_{\Gamma}$ is finitely generated and $1$-formal. 
\end{corollary}

\begin{proof}
By Corollary \ref{cor:bb mono},  
$H_{\le 2}(N_\Gamma,\Q)$ 
has trivial $\Q\Z$-action.  The conclusion 
follows from Theorem \ref{thm:ker formal}.
\end{proof}

We conclude with some examples of finitely generated Bestvina-Brady 
groups which are $1$-formal, yet admit no finite presentation.

\begin{example}
\label{ex:rp2}
Let $L=\Delta_\G$ be a flag triangulation of the real projective plane, 
$\RP^2$.  Clearly, $\Gamma$ is connected. On the other 
hand, $H_1(L,\Z)=\Z_2$, and so, by \cite{BB}, $N_{\Gamma}$ 
is not finitely presented.  But $H_1(L,\Q)=0$, and so, by 
Corollary \ref{cor:bb formal}, $N_{\Gamma}$ is $1$-formal. 
\end{example}

\begin{example}
\label{ex:poin}
Let $L=\Delta_{\Gamma}$ be a flag triangulation of a spine 
of the Poincar\'{e} homology sphere.  In \cite{BB}, Bestvina 
and Brady noted the following facts about the group $N_{\Gamma}$:
it is of type $\FP_{\infty}$ (since $\widetilde{H}_*(L,\Z)=0$), but 
not finitely presented (since $\pi_1(L)\ne 0$).  Our Corollary 
\ref{cor:bb formal} shows that $N_{\Gamma}$ is $1$-formal.  

Before finishing, we cannot but recall from \cite{BB} the 
following striking alternative about this group: either 
$N_{\Gamma}$ is a counterexample to the Eilenberg--Ganea 
conjecture, or the Whitehead conjecture is false. It would be 
interesting to know whether the formality property of $N_{\Gamma}$ 
can play a role in deciding the Bestvina--Brady alternative. 
\end{example}

\begin{ack}
This paper was started while the two authors visited the 
Mathematical Sciences Research Institute in Berkeley, 
California, in Fall, 2004, and the Abdus Salam International 
Centre for Theoretical Physics in Trieste, Italy, in Fall, 2006.  
We thank both institutions for their support and excellent facilities.

A substantial portion of the work was done during the second 
author's visit at the Institute of Mathematics of the Romanian 
Academy in October 2007.  He thanks the Institute for its 
support and hospitality during his stay in Bucharest, Romania.
\end{ack}

\vspace{-2pc}
\newcommand{\arxiv}[1]
{\texttt{\href{http://arxiv.org/abs/#1}{arxiv:#1}}}
\renewcommand{\MR}[1]
{\href{http://www.ams.org/mathscinet-getitem?mr=#1}{MR#1}}

\end{document}